\documentclass[final]{siamltex}

\usepackage{graphicx}
\usepackage{color}
\usepackage{amsmath}
\usepackage{amssymb}
\usepackage{amstext}
\usepackage[english]{babel}
\usepackage{subfig}
\newtheorem{remark}[theorem]{Remark}

\graphicspath{{./}}

\newcommand{\nph}{_{n+\frac{1}{2}}}

\begin{document}

\title{Single to Double Mill Small Noise Transition via Semi-Lagrangian Finite Volume Methods}

\author{J. A. Carrillo \footnotemark[1] \and A. Klar\footnotemark[2]
  \footnotemark[3] \and A. Roth \footnotemark[2]  }

\footnotetext[1]{Department of Mathematics, Imperial College London,
South Kensington Campus, London SW7 2AZ, UK (carrillo@imperial.ac.uk)}  
\footnotetext[2]{Technische Universit\"at Kaiserslautern, Department of Mathematics, Erwin-Schr\"odinger-Stra{\ss}e, 67663 Kaiserslautern, Germany 
  (\{klar,roth\}@mathematik.uni-kl.de)}
\footnotetext[3]{Fraunhofer ITWM, Fraunhoferplatz 1, 67663 Kaiserslautern, Germany} 
 
 \maketitle
 
\begin{abstract}
We show that double mills are more stable than single mills under stochastic perturbations in swarming dynamic models with basic attraction-repulsion mechanisms. In order to analyse accurately this fact, we
will present a numerical technique for solving kinetic mean field equations for swarming dynamics. Numerical solutions of these equations for different sets of parameters will be presented and compared to microscopic and macroscopic results. As a consequence, we numerically observe a phase transition diagram in term of the stochastic noise going from single to double mill for small stochasticity fading gradually to disordered states when the noise strength gets larger. This bifurcation diagram at the inhomogeneous kinetic level is shown by carefully computing the distribution function in velocity space.
\end{abstract}

\section{Introduction}
In the last decade, the theoretical and mathematical biology communities have paid a great deal of attention to explain large scale structures in animal groups. Coherent structures appearing from seemingly direct interactions between individuals have been reported in many different species of fish, birds, and insects \cite{HW,parrish, mogilner,BDT,couzin,camazine} and many others, see the reviews \cite{review,review2,kolo} for more literature in the subject. There are two type of patterns appearing regularly for different species called flocking and milling behavior. In the flocking behavior, individuals agree in moving in certain direction with some spatial shape that changes in time due to the effect of wind, hydrodynamic advantage, own desire, predator presence, or roosting behavior. In the milling case, individuals organise themselves rotating around certain location, that may move in time, forming a thick group in an annular/ring region. 

Many individual based (particle) models have been proposed to explain or reproduce these patterns. In general, swarming dynamics have been described by systems of interacting ordinary or stochastic differential equations: see \cite{HW,chate,couzin} for the biological aspects of the modeling of swarm dynamics and  \cite{LR,DCBC,BTTYB,BEBSVPSS,KH,HK,rome,HCH,LLE,LLE2,CKJRF} for including effects like the  interaction between individuals, self propulsion, roosting, and stochastic forces in these equations. All of these models start by including three basic mechanisms: attraction, repulsion, and reorientation. The form in which each of these effects is including in the modelling depends on the species and particular biological parameters. 

Dealing with large particle systems is cumbersome and even prohibitive numerically if the system is very large. Therefore, many authors have proposed to use a mean field approach in which one can derive effective PDEs of Fokker-Planck or Vlasov type \cite{DM1,HT08,HL08,CDP,CFRT,CCR,hauptpaper,CDP} from these microscopic systems of equations. As usual in kinetic theory, scaling arguments on these kinetic PDEs lead to macroscopic approximations of the mean field equations as in \cite{CDP,HT08,hauptpaper,MT,DFL}. 

A very interesting question that arises naturally is the stability of these patterns under perturbations. This question has been recently analysed in \cite{bigring,ABCV,CHM2} both for flocks and rings in terms of initial data for the model introduced in \cite{DCBC}. However, an stochastic perturbation also leads to very interesting phenomena called phase transition. This phase transition has been reported for the first time in swarming dynamics in the so-called Vicsek model \cite{VCBCS}, and it has been studied in very detailed at the macroscopic level in \cite{DFL}. Essentially, in these works they analyse how the system undergoes a transition from an ordered state (pattern) to a disordered state (chaos) by increasing the noise strength in the system.

The aim of our work is twofold. On the one hand, we focus on a very accurate numerical solution technique for the kinetic mean field equations in swarming dynamics. We will describe a splitting scheme with a Semi-Lagrangian solver in space and a semi-implicit finite-volume scheme in velocity space. On the other hand, we make use of this accurate numerical tool to show for the first time, up to our knowledge, a phase transition for an inhomogeneous mean field equations in swarming dynamics. By varying the amplitude of the stochastic forces, we show the influence on milling patterns in the model without noise introduced in \cite{DCBC}. The results of the microscopic, kinetic mean-field, and macroscopic equations are compared and discussed thoroughly. The numerical scheme is an improvement of the schemes used in \cite{numerikfaden,semilagrange} and of the FVM solver described in \cite{evgeniy,evgeniyhauptpaper}. As a conclusion of our study, we have discovered that double mills are quite robust and stable for small stochasticity. Mill solutions immediately leads to double mills solutions under small noise, then fading gradually toward disordered states for large noise strength. This single to double mill noise induced phase transition is the main theoretical biology implication of this work.

The paper is organized as follows. In Section \ref{models}, we quickly summarize the microscopic, kinetic mean field, and macroscopic equations derived for the model introduced in \cite{DCBC}. It contains the discussion of different regimes and associated macroscopic approximations ranging from a situation with zero stochastic force, to an intermediate case, where stochastic and propulsion force balance each other, leading finally to the dominating force case. In Section \ref{Numerical scheme}, the numerical scheme is described in details and an investigation on the order of convergence of the scheme is presented both for the homogeneous and inhomogeneous cases. Finally, Section \ref{Numerics} show the phase transition of single to double mills for small noise strength and the phase transtion toward disordered state for larger stochastic force. We first conduct a detailed comparison of microscopic and kinetic simulations for both single and double mills with and without noise for showing the accurate comparison of the models.


\section{Model Hierarchy: Microscopic, Mean-Field, and Macroscopic Models} 
\label{models}
We  consider classical models for swarming dynamics, which include forces resulting from social interaction between individuals, self-propulsion, and friction as in \cite{DCBC}. The combined effect of self-propulsion and friction is to impose an asymptotic speed for individuals. This means that individuals travel at a typical cruise speed asymptotically, similar to other classical models \cite{VCBCS,DFL}.
Furthermore, a stochastic term given by white noise is included as in \cite{CDP,review} to account for random and small error interactions. Moreover, a term describing roosting behaviour is  added as discussed and introduced in \cite{hauptpaper}. The microscopic equations are given by
\begin{align}
\label{micro}
dx_i=&\,v_i\;dt\\
dv_i=&\,v_i(\alpha-\beta|v_i|^2)\;dt 
-\frac{1}{N}\sum_{i\neq j}\nabla_{x_i} U(x_i-x_j)\;dt \nonumber\\
&- (\vert v_i\vert^2 I- v_i \otimes v_i) \nabla_{x_i} \phi (x_i)\;dt 
-\frac{A^2}{2}v_i\;dt+A\;dW_{i,t}, \label{micro2}
\end{align}
where $x_i,v_i\in\mathbb{R}^n,  i=1 \cdots N, n =2,3$. $U(x)=U(|x|)$ is a pairwise radial interaction potential. A classical example is given by  the Morse potential
$U(r)=-C_a\exp(-r/l_a)+C_r\exp(-r/l_r)$, with attraction/repulsion strengths $C_a$/$C_r$ and radius of interaction $l_a$/$l_r$. Other examples for  interaction potentials are considered in \cite{Carrillo2013,CHM} and the references therein. $\alpha$ and $\beta$ are the self propulsion parameter and $\phi$ is a potential defining the roosting force. The parameter $A>0$ describes the amplitude of the stochastic force. The main purposes of the paper are to propose an accurate numerical deterministic scheme for the mean field approximation of this microscopic system and  to study numerically the influence of the noise parameter $A$ on particular patterns present in this model. In fact, the microscopic model without noise has been shown to be paradigmatic for self-organization, since it exhibits rich dynamical patterns such as flocks, single and double mills, see \cite{DCBC,CDMBC,CDP,hauptpaper}. The main issue in the present paper is to analyse the qualitative change of the pattern as noise increases. One may expect some of these collective motions to survive for small noise up to some critical value for which the introduced stochasticity destroys the self-organized pattern. 

A classical derivation procedure as described, for example, in \cite{CDP,BH,neunzert,herleitungallgemein2,herleitungallgemein,dobru,spohn2,BCC} yields the mean field equations. For the one-particle distribution function $f(x,v,t) $ and  the density 
 $$
 \rho(x,t) = \int f(x,v,t) dv 
 $$ 
with the normalization condition $\int \rho(x,t) dx =1$, one obtains
\begin{align}
\partial_t f+\nabla_x\cdot(v f)+Sf=L f \label{fpstandard}
\end{align}
where
\begin{align}
\label{force}
Sf&=\nabla_v\cdot\big(v(\alpha-\beta|v|^2)f-(\nabla_x U*\rho)f- (\vert v \vert^2 I - v \otimes v )\nabla_x \phi f\big)
\end{align}
and
\begin{align}
\label{stochasticforce}
Lf&=\frac{A^2}{2}\nabla_v\cdot (vf+\nabla_vf).
\end{align}

Different regimes and approximate solutions for macroscopic quantities can be obtained from the above mean field equation using different scaling assumptions ranging from very small to very large stochastic force $A$.

\subsection{Vanishing stochastic force}
\label{A=0}
The deterministic case $A=0$ has been treated in detail in many publications. A special feature is the appearance of single and double mill solutions, we refer to \cite{CDP}. Single mills of the mean field equation 
\begin{align*}
\partial_tf+\nabla_x\cdot(v f)+Sf=0, 
\end{align*}
where $S$ is given by \eqref{force}, can be found as follows.
We look for a  mono-kinetic solution, compare \cite{CDP,CDMBC}
$$
f (x,v) = \rho(x) \delta_{u(x)} (v).
$$
$\rho$ and $u$ are found by integrating against $dv $ and $v \; dv$ and closing the equations with the above Ansatz. One obtains
\begin{eqnarray*}
  \partial_t \rho   +   \nabla_{x} \cdot  (\rho u)  = 0
\end{eqnarray*}
and the momentum equation
\begin{align*}
\partial_t u + (u \cdot  \nabla_{x}  ) u = u ( \alpha
- \beta  \vert u\vert^2 ) - \nabla_x U \star \rho - 
(\vert u \vert^2 I -   u  \otimes u)  \nabla_x \phi 
\end{align*}
on the support of the density $\rho$. Stationary distributions can be found in the following way. Assuming $\beta \vert u \vert^2 = \alpha$, we obtain from the hydrodynamic equations above that
\begin{eqnarray*}
   \nabla_{x} \cdot  (\rho u)
  &=& 0\\
  (u \cdot  \nabla_{x}  ) u
 & =&
  - \nabla_x U \star \rho -  (\vert u \vert^2 I -   u  \otimes u)  \nabla_x \phi  
\end{eqnarray*}
on the support of the density $\rho$. Assuming a rotatory solution
given by
$$
u = \sqrt{\frac{\alpha}{\beta}}   \frac{x^{\perp}}{\vert x \vert},
$$
compare \cite{CDP}, and looking for radial densities $\rho =
\rho(\vert x \vert) $, we obtain that the continuity equation is trivially satisfied and
since
\begin{eqnarray*}
  (u \cdot  \nabla_{x}  ) u
  = - \frac{\alpha}{\beta}  \frac{x}{\vert x \vert^2}
\end{eqnarray*}
then
\begin{eqnarray*}
 - \frac{\alpha}{\beta}  \frac{x}{\vert x \vert^2}    = - \nabla_x U \star \rho  - \frac{\alpha}{\beta} \frac{x}{\vert x \vert^2}  x \cdot  \nabla_x \phi. 
\end{eqnarray*}
Assuming $\phi(x) = \phi(\vert x \vert)$ we end up with an
integral equation for $\rho$:
\begin{eqnarray}
\label{integral}
 U \star \rho  = D   + \frac{\alpha}{\beta}  \left(\ln \vert x \vert -     \phi (\vert x \vert)\right)
\end{eqnarray}
in the support of the density $\rho$. A numerical investigation of these stationary states of the hydrodynamic equations, called single mills, is performed in \cite{hauptpaper}.

A so called double mill solution can also be derived from an appropriate Ansatz 
with two opposing kinetic velocities, see \cite{CDP}. In short, given the spatial profile of a single mill, there is always a double mill solution with the same density but changing the sign of the velocity to part of it. Let us finally remark that in particle simulations of the microscopic system \eqref{micro}-\eqref{micro2}, double mills are observed very often by peforming large perturbations of single mills when $A=0$. Moreover, these double mills are typically interlaced instead of overlapping on the same spatial region. This means that the velocity of one mill is not totally opposite of the velocity of the other mill at the same spatial point.

\subsection{Balance of stochastic and propulsion forces}
Assume that $A$ and $\alpha,\beta$ are large for long time scales. The corresponding scaling yields the following equations
\begin{align*}
\epsilon \partial_tf+\nabla_x\cdot(v f)+ \tilde{S}f=  \frac{1}{\epsilon} \tilde{L}f 
\end{align*}
where
\begin{align*}
\tilde{L}f&=-\nabla_v\cdot\big(v(\alpha-\beta|v|^2)f \Big)  + \frac{A^2}{2}\nabla_v\cdot (vf+\nabla_vf)
\end{align*}
and
\begin{align*}
\tilde{S}f&=\nabla_v\cdot\big(-(\nabla_x U*\rho)f- (\vert v \vert^2 I - v \otimes v )\nabla_x \phi f\big).
\end{align*}
To zeroth order in an expansion in the small parameter $\epsilon$, we have $\tilde{L}(f_0) =0$. The kernel of $\tilde{L}$ is given by the solution of 
\begin{align*}
\nabla_v\cdot(v(\alpha-\beta |v|^2)f)=\frac{A^2}{2}\nabla_v\cdot(vf+\nabla_vf).
\end{align*}
A simple computation gives that the local equilibrium solution reads as 
\begin{align*}
f_0(x,v)=\rho_0 C \exp\left(-\frac{2}{A^2}\psi(v)\right) \quad \mbox{with } \psi(v)=\frac{2}{A^2}\left(\beta\frac{|v|^4}{4}- (\alpha - \frac{A^2}{2})\frac{|v|^2}{2}\right)
\end{align*}
with a normalizing constant $C$ and a given spatial density $\rho_0(x,t)$.

\

\begin{remark}
For $A$ going to zero $f_0$ goes to 
$$\rho_0 \delta_{\vert v \vert^2 = \alpha/\beta} (v).$$
For $A$ going to infinity $f_0$ converges to 
a standard Maxwellian $ \rho_0 M(v) $.
\end{remark}

\

To compute $\rho_0$ one has to proceed to first order. One obtains 
\begin{align*}
 \tilde{L} f_1  = g = v \nabla_x f_0 - \mbox{div}_v ((\nabla_x U \star \rho_0) f_0 ) -  \mbox{div}_v ((\vert v \vert^2 I - v \otimes v )\nabla_x \phi f_0).
\end{align*}
This equation is solvable since  
$
\int g dv =0.
$
Then, it is not difficult to check that
\begin{align*}
-\int v f_1  dv = D_0  \nabla_x \rho_0 + D_1 (\nabla_x U \star \rho_0) \rho_0  + D_2 \nabla_x \phi \rho_0
\end{align*}
with the matrices $D_i$ given by
\begin{align*}
D_i =  \int v \otimes \tilde{L}^{-1}( g_i(v)) dv,\; i=0,1,2
\end{align*}
and the vector fields
\begin{align*}
g_0 &=  -v C \exp\left(-\frac{2}{A^2} \psi(v)\right)\,,\\
g_1 &= -C \nabla_v  \left[\exp\left(-\frac{2}{A^2} \psi(v)\right)\right]\,,\\
g_2 &= -C \nabla_v \cdot \left[ (\vert v \vert^2 I - v \otimes v) \exp\left(-\frac{2}{A^2} \psi(v)\right)\right].
\end{align*}
Integrating the scaled equation and using the above  gives
\begin{align}\label{macrob}
 \partial_t \rho_0 = \nabla_x \cdot((D_1 \nabla_x U \star \rho_0    + D_2 \nabla_x \phi) \rho_0 + D_0 \nabla_x \rho_0) .
\end{align}

\subsection{Dominating stochastic force}
\label{Alarge}
We refer to \cite{CDP,hauptpaper,BGKMW07,GKMW07,HKMO09} and consider a scaling with large $A$ and long time scales. This yields
\begin{align*}
\epsilon \partial_tf+\nabla_x\cdot(v f)+Sf=\frac{1}{ \epsilon }L f  
\end{align*}
A straightforward asymptotic expansion gives to zeroth order
$f_0 = \rho_0 M(v)$ where $M$ is the standard Maxwellian and 
\begin{align*}
 \partial_t \rho_0 = \nabla_x\cdot \left( (\nabla_x U \star \rho_0 + \nabla_x \phi) \rho + \nabla_x \rho_0\right) . 
\end{align*}
The solution of the stationary problem is given by
\begin{align}
\label{statbigA}
 \rho_0 =   \tilde{C} \exp (-  U \star \rho_0 - \phi)
\end{align}
with $\tilde{C}$ prescribing the total mass. The associated momentum is 
$$
 -\int v f_1 dv =
  \nabla_x \rho_0 + (\nabla_x U \star \rho_0) \rho_0  +\nabla_x \phi \rho_0. 
$$

To summarize this overview, we expect the stationary solutions of the mean field equations to change as follows when $A$ goes from $0$ to infinity.
For $A=0$ and well prepared inital conditions, one obtains a single mill solution, i.e. a $\delta$ solution in velocity and a radial solution in $x$ of equation  (\ref{integral}). Double mills appear depending on the initial conditions and typically for large perturbations of the mill solution. As $A$ becomes very large, one approaches a  Maxwellian solution in velocity and a solution of equation (\ref{statbigA})  in $x$. These two behaviors will be confirmed by the numerical investigations in Section 4. 

On the other hand, the transition from the mill ordered state to the disordered state by increasing the noise is not clear at this point. The macroscopic equations  in the middle regime seem to indicate a density profile near a non-isotropic maxwellian due to the matrices $D_i$ in \eqref{macrob}. However, the change in the spatial and velocity profile are not clear for small noise when approaching the monokinetic solution. We wil investigate this issue in detail in  Section 4 showing that double mills play an important role in this transition.

\section{Numerical Scheme}
\label{Numerical scheme}
\subsection{Discretisation}
In the following, we present a second order  Semi-Lagrange Finite Volume method for the Fokker-Planck equations. We describe the method in 2-D, see \cite{KRS07} for similar situations in fibre dynamics. An extension to 3-D is straightforward, but computationally  expensive, see \cite{3dfadenmodell} for similar problems. We split (\ref{fpstandard}) using  a Strang splitting  \cite{Splittingstrang}, that means we consider
\begin{align}
\partial_tf^*&=-\frac{1}{2}Sf^*+\frac{1}{2}Lf^* & f^*(t)&=f(t)\label{velpart}\\
\partial_tf^{**}&=-v\nabla_xf^{**} & f^{**}(t)&=f^*(t+\tau)\label{spatialpart}\\
\partial_tf&=-\frac{1}{2}Sf+\frac{1}{2}Lf & f(t)&=f^{**}(t+\tau) \nonumber 
\end{align}
where $S$ and $L$ are given by \eqref{force} and \eqref{stochasticforce} respectively.
The system is discretised on the grid $(x_i,y_j,v^k)$, where $x_{ij}=(x_i,y_j)\in\mathbb{R}^2$ for the spatial domain and $v^k\in\mathbb{R}^2$ are points in the velocity domain. Furthermore, we denote by $t_n$ the discretisation points in time with constant step size $\tau$. For (\ref{spatialpart}) we use a second order Semi-Lagrange method (see \cite{sem6,semilagrange,qiu}), where we apply an interpolation procedure using cubic Bezier curves. Proceeding in the usual way,
 we look back in time along the characteristic curves starting in the respective grid point and obtain a grid value for the current time step by interpolating from the old grid values at the characteristic endpoint. On each cell $[x_i,x_{i+1}]\times[y_j,y_{j+1}]$ of the grid in space, the solution at a given time is given by the polynomial reconstruction
\begin{align}
B_{i+\frac{1}{2}}(t_x,t_y)&=\sum_{k,l=0}^3B_{k3}(t_x)B_{l3}(t_y)\xi_{kl}\label{bezier}
\end{align}
with the Bernstein polynomials $B_{k3}(t)=\binom{3}{k}t^k(1-t)^k, t \in [0,1]$ and the 16 control points $\xi_{kl}$. We used the notations $t_x=\frac{x-x_i}{h_x},t_y=\frac{y-y_j}{h_y}$ here. Due to the nature of the Bernstein polynomials, the interpolant never leaves the convex hull of the control points. They are chosen appropriately so that we have a third order local approximation for smooth data and a reduced order interpolant without oscillations for non-smooth data. 

 In order to obtain a cubic Bezier polynomial like (\ref{bezier}) with appropriate control values $\xi_{kl}$, one computes an interpolating Newton polynomial and then performs a basis transformation from the Newton basis to the Bernstein polynomials. We will demonstrate the procedure for a 1D problem. Assume we have the stencil $F=(f_{i-1},f_i,f_{i+1},f_{i+2})^T$. Then, for $x\in[x_i,x_{i+1}]$ and $t=t_x = \frac{x-x_i}{h_x}$, we can write  the Newton polynomial $N_{i+\frac{1}{2}}(t)$ as follows
\begin{align*}
N(t)&=\begin{pmatrix}1&(t+1)&t(t+1)&(t-1)t(t+1)\end{pmatrix}\\
b&=\begin{pmatrix}b_0&b_1&b_2&b_3\end{pmatrix}^T\\
N_{i+\frac{1}{2}}(t)&=N(t) b
\end{align*}
where we obtain the coefficients
\begin{align*}
b&=\begin{pmatrix}
N(-1)\\N(0)\\N(1)\\N(2)
\end{pmatrix}^{-1}F=N^{-1}F.
\end{align*}

\begin{figure}
\centering
\includegraphics[width=6cm]{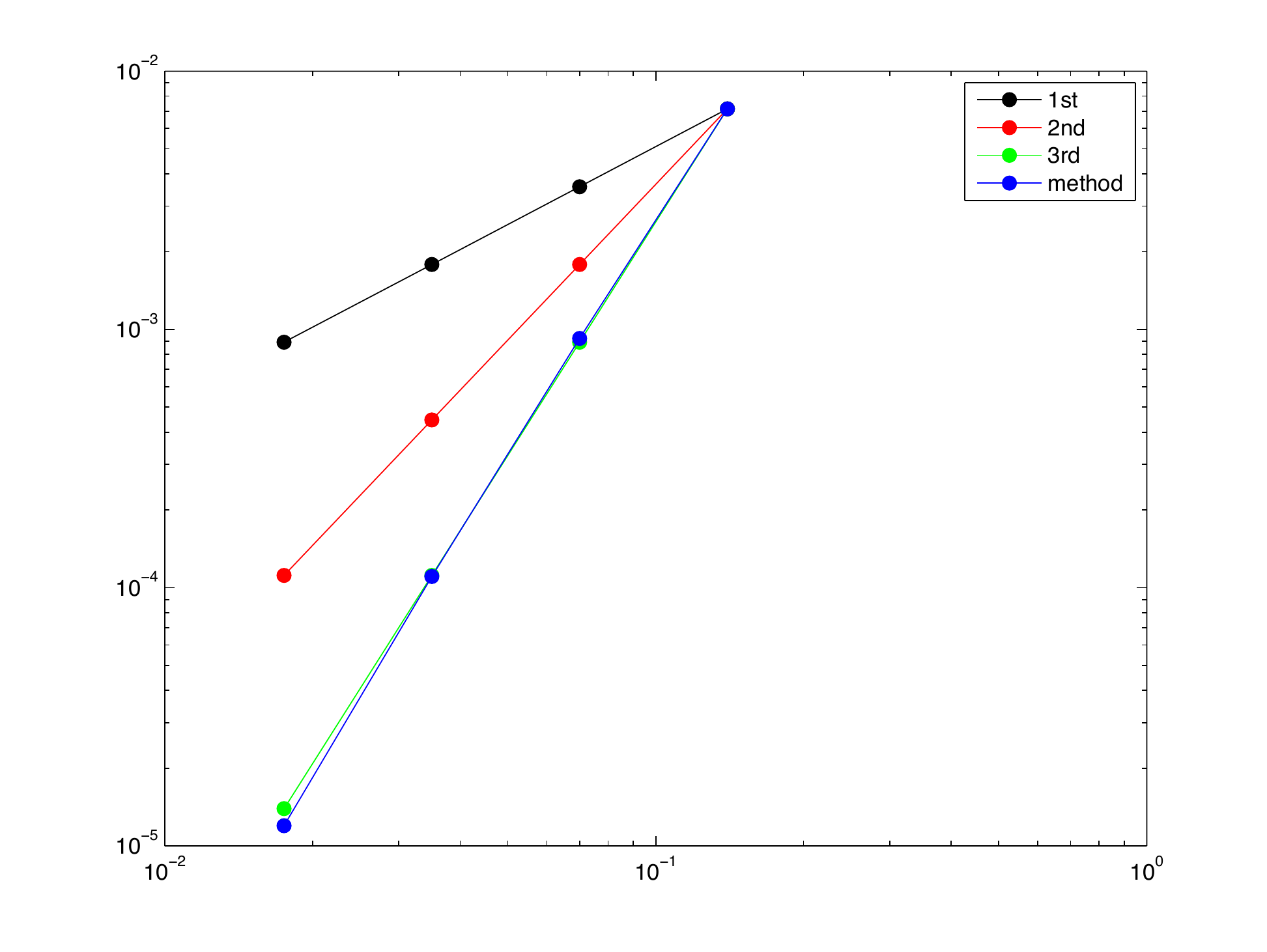}
\caption{Convergence of the Semi-Lagrange method with Bezier interpolation for linear advection as in (\ref{spatialpart}). One obtains the predicted third order method. \label{semilagconv}}
\end{figure}

\begin{figure}[Hbt]
\centering
\includegraphics[width=6.8cm]{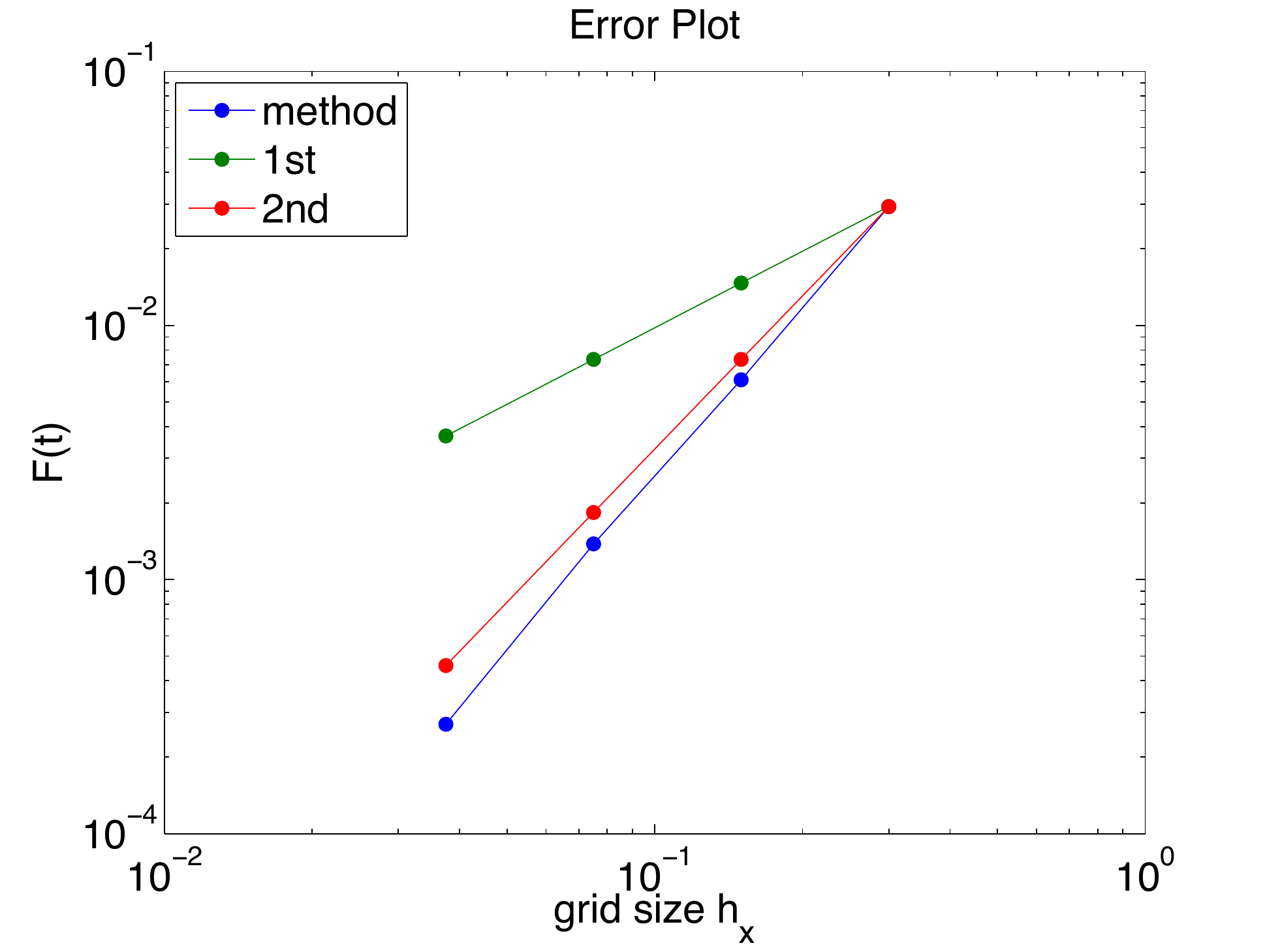}
\caption{Logarithmic error plot for different grid sizes for the finite volume method in the space homogeneous case. The coefficient $A$ was chosen as $A=0.15$. The convergence is second order.\label{konv1}}
\end{figure}

In the same way, we can write  the Standard polynomial basis, the Bernstein polynomial basis, the Bezier coefficients and the Bezier polynomial:
\begin{align*}
S(t)&=\begin{pmatrix}1&t&t^2&t^3\end{pmatrix}\\
B(t)&=\begin{pmatrix}(1-t)^3&3t(1-t)^2&3t^2(1-t)&t^3\end{pmatrix}\\
\xi&=\begin{pmatrix}\xi_0&\xi_1&\xi_2&\xi_3\end{pmatrix}^T\\
B_{i+\frac{1}{2}}(t)&=B(t)\xi
\end{align*}
Now   a change of basis gives
\begin{align*}
B(t)&=S(t)T_B \quad \mbox{ and } \quad N(t)=S(t)T_N & \Rightarrow\;\;\;\;N(t)&=B(t)T_B^{-1}T_N\,.
\end{align*}
In order to achieve the right interpolation order, we want the Bezier polynomial to be equal to the Newton polynomial, i.e. $B(t)\xi=N(t)b$. That can be achieved if we set
\begin{align*}
\xi&=T_B^{-1}T_NN^{-1}F\\
\Rightarrow\;\;\;\;\begin{pmatrix}\xi_0\\\xi_1\\\xi_2\\\xi_3\end{pmatrix}
&=\begin{pmatrix}
                             												 f_i\\
     \frac{5f_i}{6} - \frac{f_{i-1}}{9} + \frac{f_{i+1}}{3} - \frac{f_{i+2}}{18}\\
 \frac{5f_{i+1}}{6} - \frac{f_{i+2}}{9} +     \frac{f_i}{3} - \frac{f_{i-1}}{18}\\
                             											 f_{i+1}\\\end{pmatrix}
\end{align*}
Since the Bezier polynomial with these control values is the same as the Newton polynomial, we know that we have a local approximation error of fourth order for smooth data, resulting in a third order consistency error for linear advection in the Semi-Lagrange scheme, see figure \ref{semilagconv}. However, at discontinuities, oscillations will occur. Avoiding them is straight forward, due to the convex hull property. If one of the $\xi_i$ leaves the interval $[f_i,f_{i+1}]$ of neighbouring grid values, we put it back to either of the boundaries $f_i$ or $f_{i+1}$. We note that the mass conservation is destroyed by such a procedure. Thus, additionally a mass conservation procedure as in \cite{semilagrange,douglas00} has been implemented. In literature, similar examples of Semi Lagrange schemes and splitting methods for the Vlasov equation are known, see for example \cite{CV,filbet,qiu,sem6,VLR}. The contributions mostly differ in the interpolation scheme and the handling of oscillatory solutions. For the advection problem (\ref{spatialpart}), our procedure can be easily extended to 2D  by extending  the above described vectors in an appropriate way.
\\
\begin{figure}
\centering
\includegraphics[width=7cm]{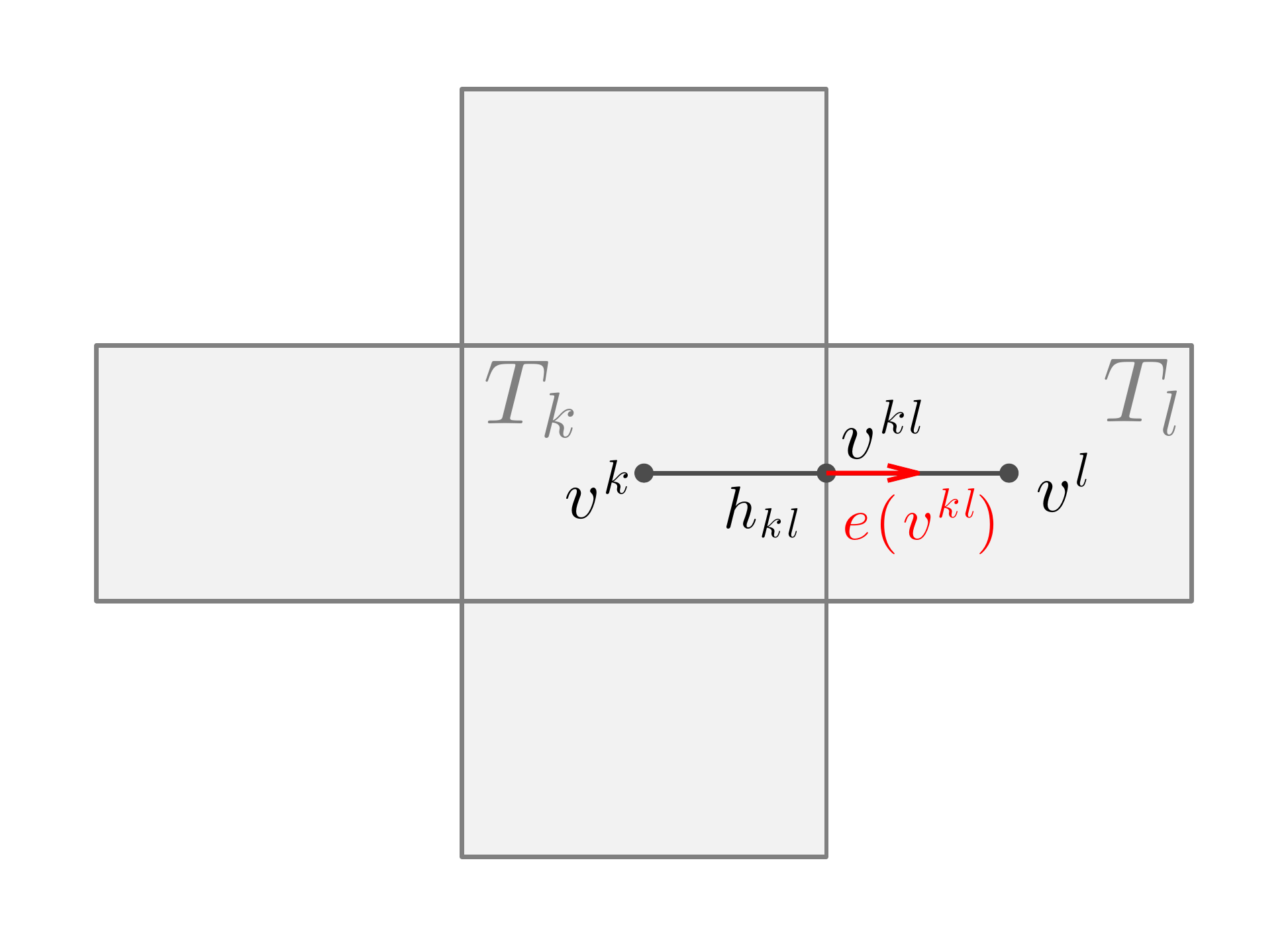}
\caption{A patch of the rectangular grid for FVM.\label{fvmvis}}
\end{figure}
For the discretization in velocity, we apply a finite volume scheme \cite{leveque} of second order with a Lax Wendroff type approximation \cite{quarteroni} for the advection, which is an enhancement of the methods used in \cite{numerikfaden,evgeniyhauptpaper}: Denote the midpoints of rectangular cells $T_k$ for (\ref{fpstandard}) by $v^k$. The boundary between cells $k$ and $l$ is called $T_{kl}$ with the normal $e(v)$ and the boundary midpoint $v^{kl}$. We further denote the distance between cell midpoints $v^k$ and $v^l$ by $h_{kl}$, which decomposes into $h_{kl}=h_1+h_2$, where $h_1,h_2$ are the distances to the boundary midpoint $v^{kl}$. Since we have an equidistant, rectangular grid in $x$ and $y$-direction, $h_{kl}$ is either $h_x$ or $h_y$ and $h_1,h_2$ are $h_x/2,h_y/2$ respectively, see Figure \ref{fvmvis}. $N(k)$ is the set of the indices of neighbour cells to cell $k$. For the FVM nethod, one computes the cell averages over the cells $[t_n,t_{n+1}]\times T_k$ by integrating (\ref{velpart}). This gives
\begin{align*}
\frac{1}{\tau}\frac{1}{|T_k|}\int_{T_k}\int_{t_n}^{t_{n+1}} \partial_tf\;dtdv=&-\frac{1}{\tau}\frac{1}{|T_k|}\int_{t_n}^{t_{n+1}}\int_{T_k}\nabla_v\cdot(F(x,v) f(v,t))\;dvdt\\
&+\frac{1}{\tau}\frac{A^2}{2|T_k|}\int_{t_n}^{t_{n+1}}\int_{T_k}\nabla_v\cdot(\nabla_vf(v,t))\;dvdt
\end{align*}
where the notation $F(x,v):=v(\alpha-\beta|v|^2)-(\nabla_x{U}*\rho)-(\nabla_x \phi \cdot v^\perp)v^\perp-\frac{A^2}{2}v$ is used. Rewriting
\begin{align*}
f_k^n=\frac{1}{|T_k|}\int_{T_k}f(v,t_n)\;dv
\end{align*}
and applying the midpoint rule for the time integration and  the divergence theorem for the integration w.r.t $v$, one gets
\begin{align}
\frac{f_k^{n+1}-f_k^n}{\tau}=&\frac{1}{|T_k|}\sum_{l\in N(k)}\int_{T_{kl}}f(v(S),t\nph)F(x,v(S))\cdot e(v(S))\;dv(S)\nonumber\\
&+\frac{A^2}{2|T_k|}\sum_{l\in N(k)}\int_{T_{kl}}\nabla_vf(v(S),t\nph)\cdot e(v(S))\;dv(S)+\mathcal{O}(\tau^2)\nonumber\\
\begin{split}=&\frac{1}{|T_k|}\sum_{l\in N(k)}|T_{kl}|f(v^{kl},t\nph)F(x,v^{kl})\cdot e(v^{kl})\\
&+\frac{A^2}{2|T_k|}\sum_{l\in N(k)}|T_{kl}|\nabla_vf(v^{kl},t\nph)\cdot e(v^{kl})+\mathcal{O}(|T_{kl}|^3)+\mathcal{O}(\tau^2)\,.\end{split}\label{discr}
\end{align}
Now, we have to approximate the quantities $f(v^{kl},t\nph)$ and $\nabla_vf(v^{kl},t\nph)\cdot e(v^{kl})$.  If we use
\begin{align*}
f(v^{kl},t\nph)&=\frac{1}{2}\Big[(f_k^n+f_l^n)-\frac{\tau}{h_{kl}} e(v^{kl})\cdot\Big(F(x,v^l)f_l^n-F(x,v^k)f_k^n\Big)\Big]:=f^{LW}\,,
\end{align*}
we end up with a second order Lax-Wendroff flux function. As any higher order flux function, it will cause oscillations for non-smooth solutions and one has to apply a  limiter scheme to remedy this. In the present example, we will mix the Lax-Wendroff flux with a first order upwind advection $f^{UW}$ for non-smooth solutions:
\begin{align*}
f(v^{kl},t\nph)&=\phi(\theta_{kl})f^{LW}+(1-\phi(\theta_{kl}))f^{UW}
\end{align*}
where
\begin{align*}
f^{UW}&=\begin{cases}
f_k^n&\mbox{,  if }F(x,v^{kl})\cdot e(v^{kl})>0\\
f_l^n&\mbox{,  else}
\end{cases}
\end{align*}
and $\theta_{kl}$ measures the smoothness of the solution. $\phi$ is the van-Leer limiter \cite{vanleer,toro}
\begin{align*}
\phi(\theta)&=\frac{\theta+|\theta|}{1+|\theta|}.
\end{align*}
Using this limiting procedure, oscillations are avoided. In order to get a value for $\nabla_vf\cdot e(v^{kl})$, one can expand $f$ around $v^{kl}$ in edge normal direction and obtains
\begin{align*}
\nabla_vf(v^{kl})\cdot e(v^{kl})&=\frac{f(v^l)-f(v^k)}{h_{kl}}+\mathcal{O}(h_{kl}^2)
\end{align*}
for the equidistant grids used here. Together with the fluxes from the other $l\in N(k)$, the second order terms vanish and one obtains an estimate of oder  $\mathcal{O}(h_{kl}^3)$. In total, multiplying by $|T_{kl}|$ and dividing by $|T_i|$, we end up with a second order approximation on the velocity grid. We have to note however, that the discretization is not second order in time, since we did not account for the evaluation of $\nabla_vf\cdot e(v^{kl})$ at $t\nph$. This can be remedied by applying the trapezoidal rule in (\ref{discr}) instead of the midpoint rule for the time integral of the diffusion terms, which would lead to a Crank-Nicholson scheme for the diffusive part of the problem. The convergence rate in time does not influence the numerical results that much, since the overall error is still very near to a second order scheme: For a numerical confirmation of the convergence rate of the velocity discretization in the case without interaction and roosting, see Figure \ref{konv1}.

\begin{figure}[Hbt]
\centering
\includegraphics[width=6.2cm]{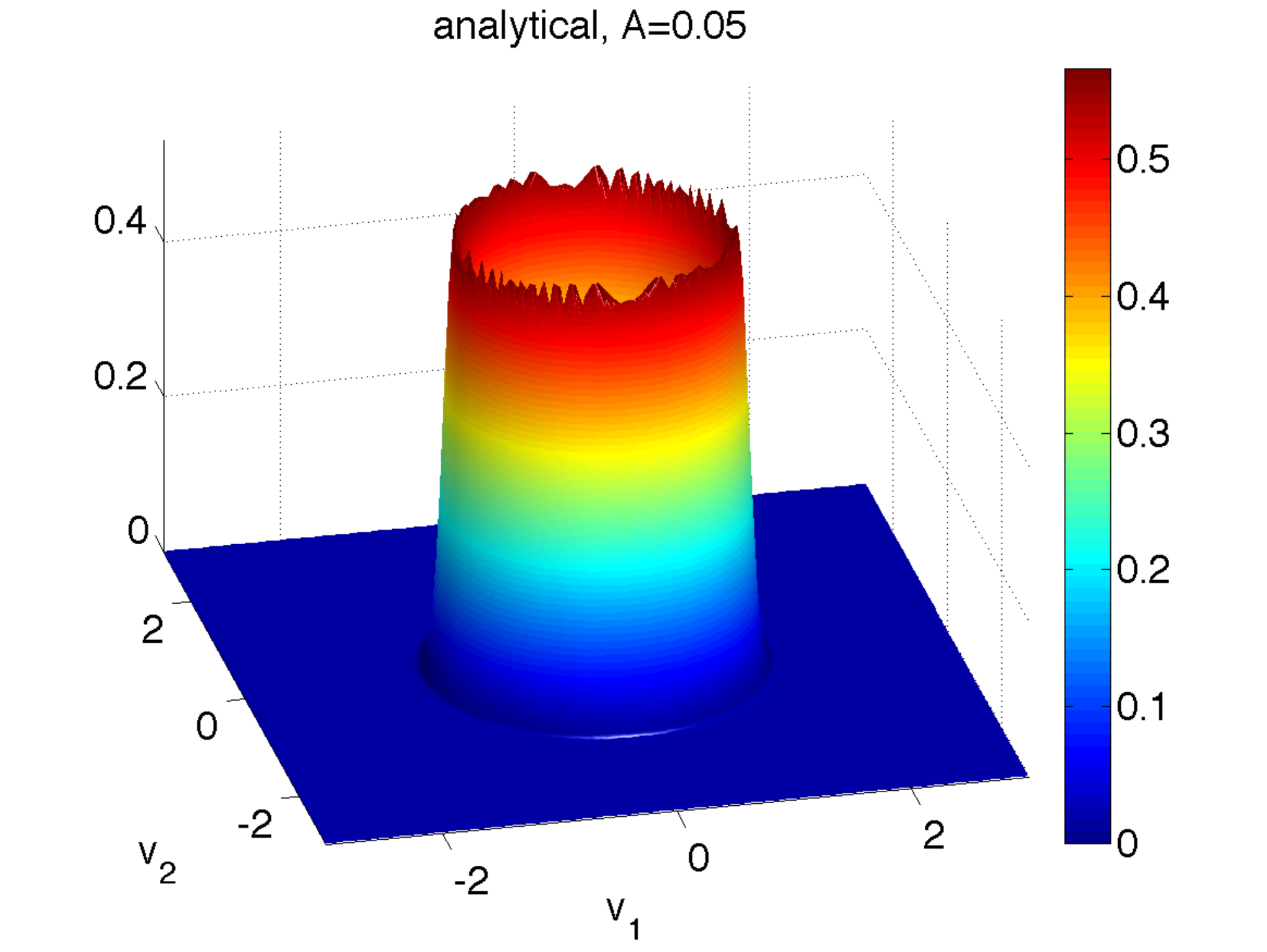}
\includegraphics[width=6.2cm]{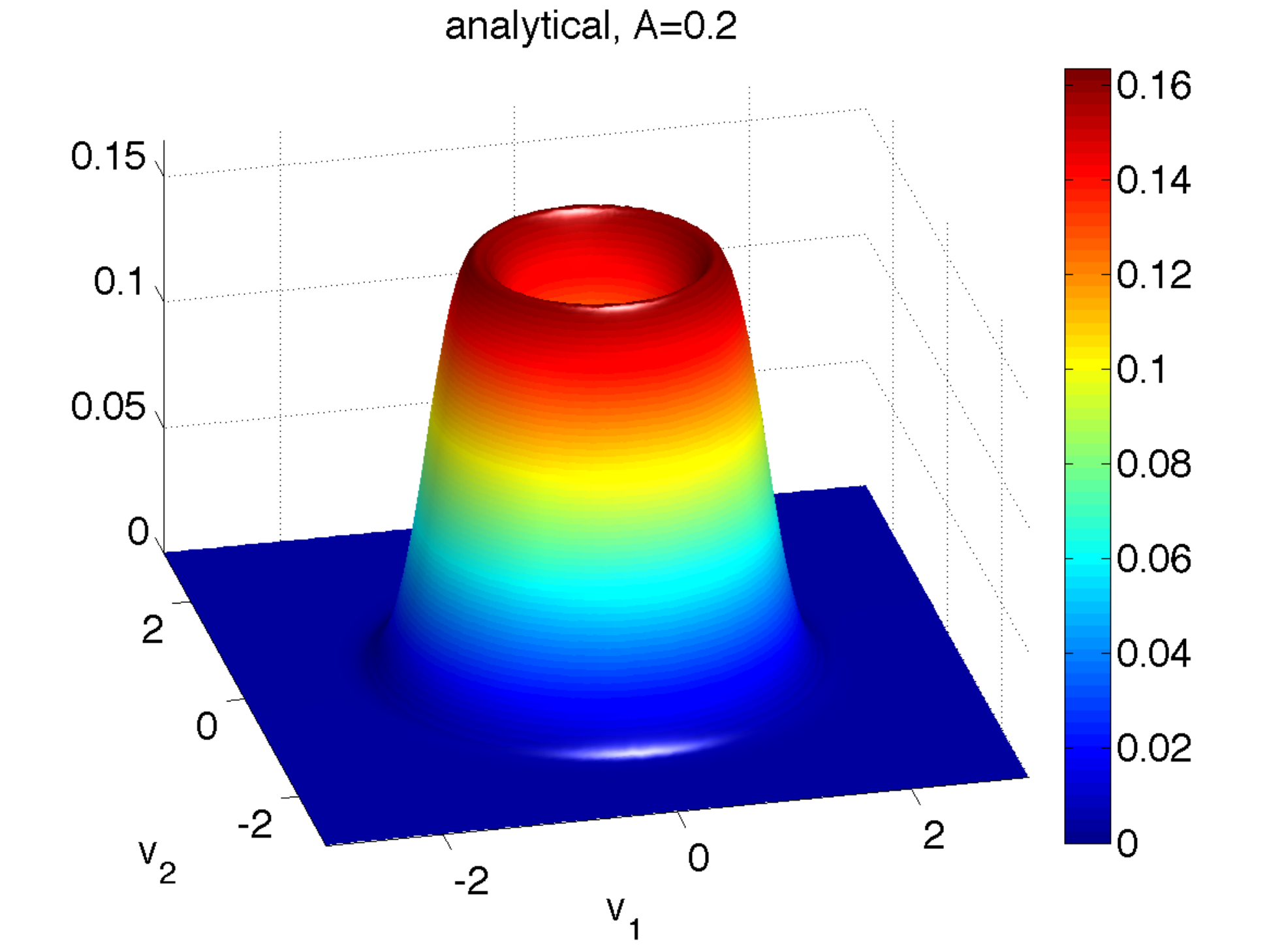}\\
\includegraphics[width=6.2cm]{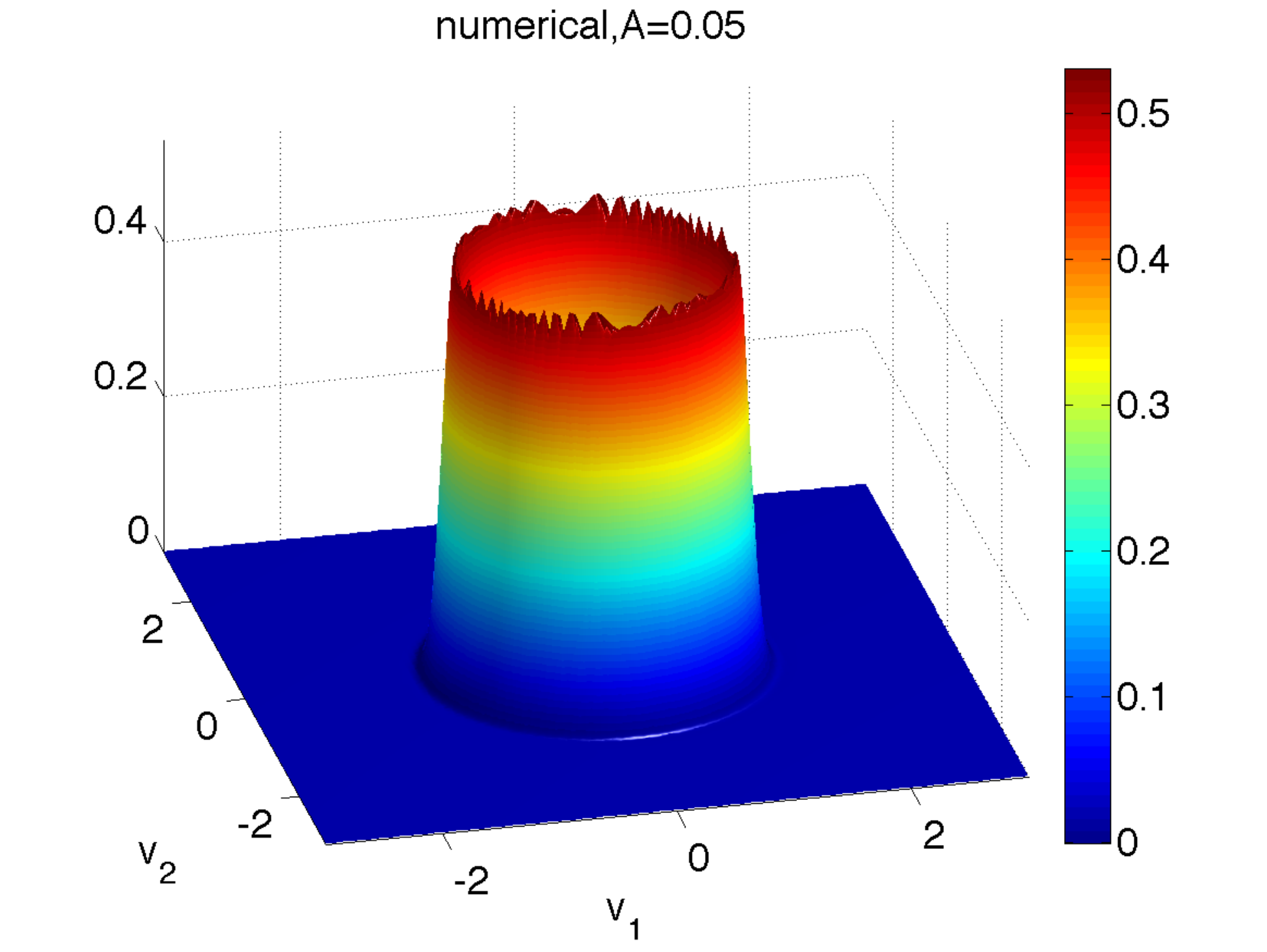}
\includegraphics[width=6.2cm]{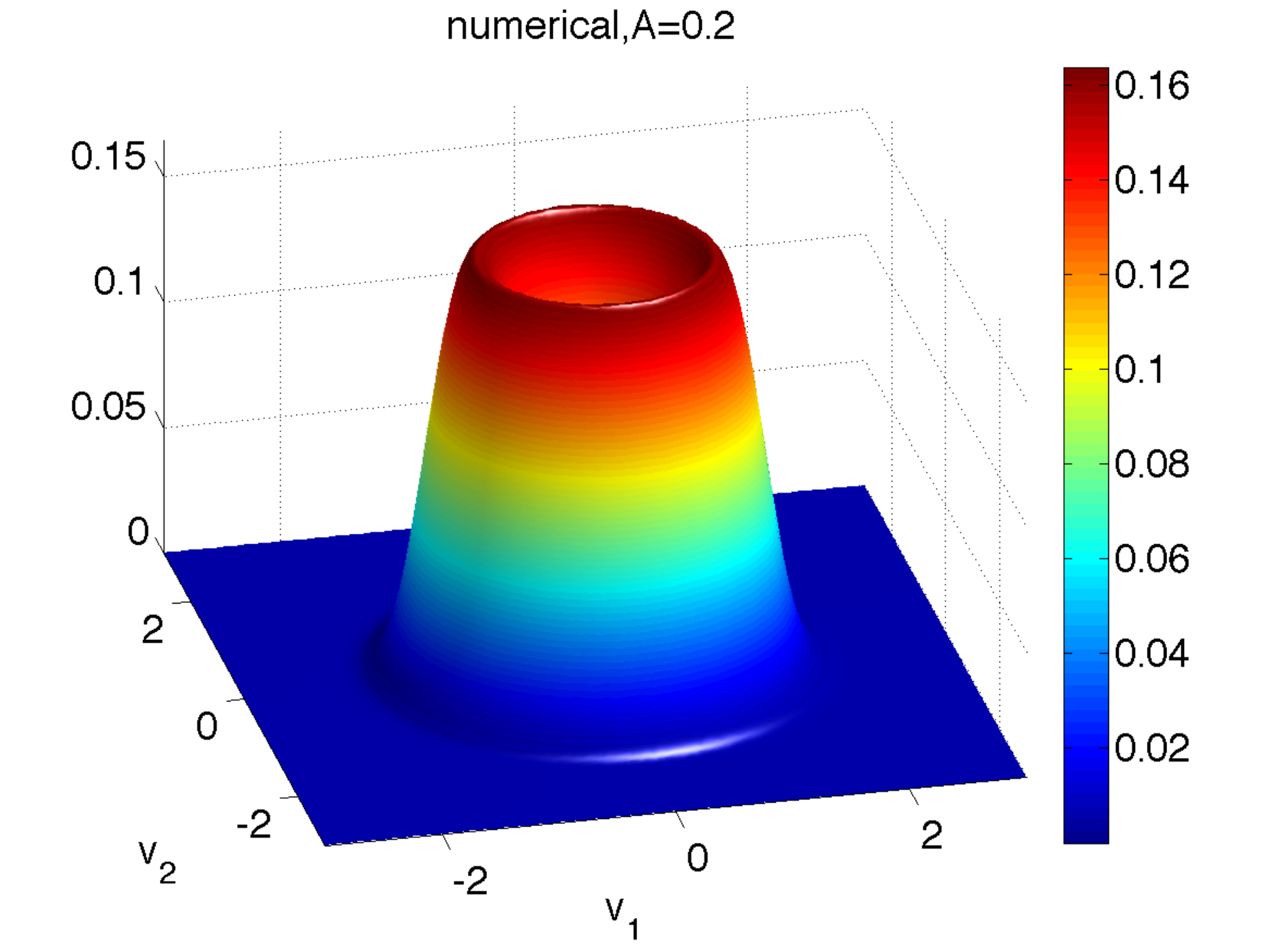}\\
\caption{On the left, the diffusion coefficient $A$ is  set to $0.05$, on the right, it is  $A=0.2$. From top to bottom we have: Analytical equilibrium $f_{eq}$ and numerical equilibrium $f_{eq}^N$. For $A\rightarrow 0$, $f_{eq}$ will approach $\delta(|v|-\sqrt{\frac{\alpha}{\beta}})$.
\label{vergleich}}
\end{figure}

In matrix form   (\ref{discr}) can be rewritten as
\begin{align}
f^{n+1}&=f^n-\tau\cdot\Big(A_Tf^n-A_Df^n\Big)\label{semidis}
\end{align}
with matrices $A_T$ for the transport coefficients and $A_D$ for the diffusion coefficients. For large values of the diffusion parameter ($A>1$) one will need a very small time step, if the above method (\ref{semidis}) is to converge. So instead of (\ref{semidis}), we choose the semi-implicit ansatz
\begin{align}
f^{n+1}&=f^n-\tau\cdot\big(A_Tf^n-A_Df^{n+1}\Big)\nonumber
\end{align}
or
\begin{align}
 (I-\tau A_D)f^{n+1}&=(I-\tau A_T)f^n.\label{linsystem}
\end{align}
$I-\tau A_D$ can be proven to be strictly diagonally dominant and has positive diagonal entries\label{sec:impldiff}. 
We use the conjugate gradient method for solving system (\ref{linsystem}). For the implementation, it is unavoidable to use some sort of sparse format for the matrices, which decreases computing time by a large factor, since we have to solve this system at every point $(x^i,y^j)$ in every time step.

\subsection{Numerical investigation of the space homogeneous case}
First we test our finite volume scheme for the space homogeneous  equation without interaction and roosting force. The equation reads 
\begin{align}
\partial_tf+\nabla_v\cdot(v(\alpha-\beta |v|^2)f)=\frac{A^2}{2}\nabla_v\cdot(vf+\nabla_vf).\label{nurgeschw}
\end{align}
For the computations we use the numerical 
 values  $\alpha= 0.07$ and $ \beta = 0.05$.
The analytical equilibrium solution is given by 
\begin{align*}
f_{eq}(v)=C\cdot\exp\left(-\left(\tilde{\beta}\frac{|v|^4}{4}-\tilde{\alpha}\frac{|v|^2}{2}\right)\right)
\end{align*}
with
$\tilde{\beta}=2 \beta/A^2$ and $\tilde{\alpha}=2 \alpha /A^2 -1$. For  $A\neq 0$ the solution of the time dependent Fokker-Planck equation (\ref{nurgeschw}) converges to this equilibrium state, compare figure \ref{vergleich}. The convergence to equilibrium is investigated looking at the
distance between numerical $f_{eq}^{N}$ and analytical equilibrium  $f_{eq}$ in $L^2$.  The values of
this functional
for different grid sizes are displayed in Figure \ref{konv1}. 

\subsection{Numerical convergence analysis for the full problem}
Up to now, we only tested the numerical scheme for spatial and velocity domain separately. In this section  a numerical convergence analysis is performed for the two regimes of large diffusion $A=3.0$ and no diffusion $A=0.0$, for the full problem in space and velocity. Starting from a Single Mill initial condition in all cases, we will compare  the stationary distributions for different grid sizes.  Let $\rho=\int f \;dx$ be the density of the numerically computed stationary distribution function and  $\rho_R$ be the reference density. To measure the errors we compute the distance 
\begin{align}
|\rho-\rho_R|_2&=\sqrt{\int(\rho-\rho_R)^2\;dx}.\label{norms}
\end{align}
The reference solutions are given on a finer grid than the numerical solutions. We will interpolate the numerical data with a 3rd order accurate procedure to the reference grid and compute the distance there. The occurring integrals will then be evaluated with the midpoint rule on each grid cell.

\subsubsection{A=3.0}
In this case, we do not have an analytical equilibrium. However, a very good approximation of the stationary  density $\rho_R$ can be obtained by taking  the solution of the fixed-point equation
\begin{align*}
\rho_R&=\frac{\exp(-U*\rho_R)}{\int \rho_R dx}
\end{align*}
 
\begin{figure}[hbt]
\centering
\includegraphics[width=\textwidth]{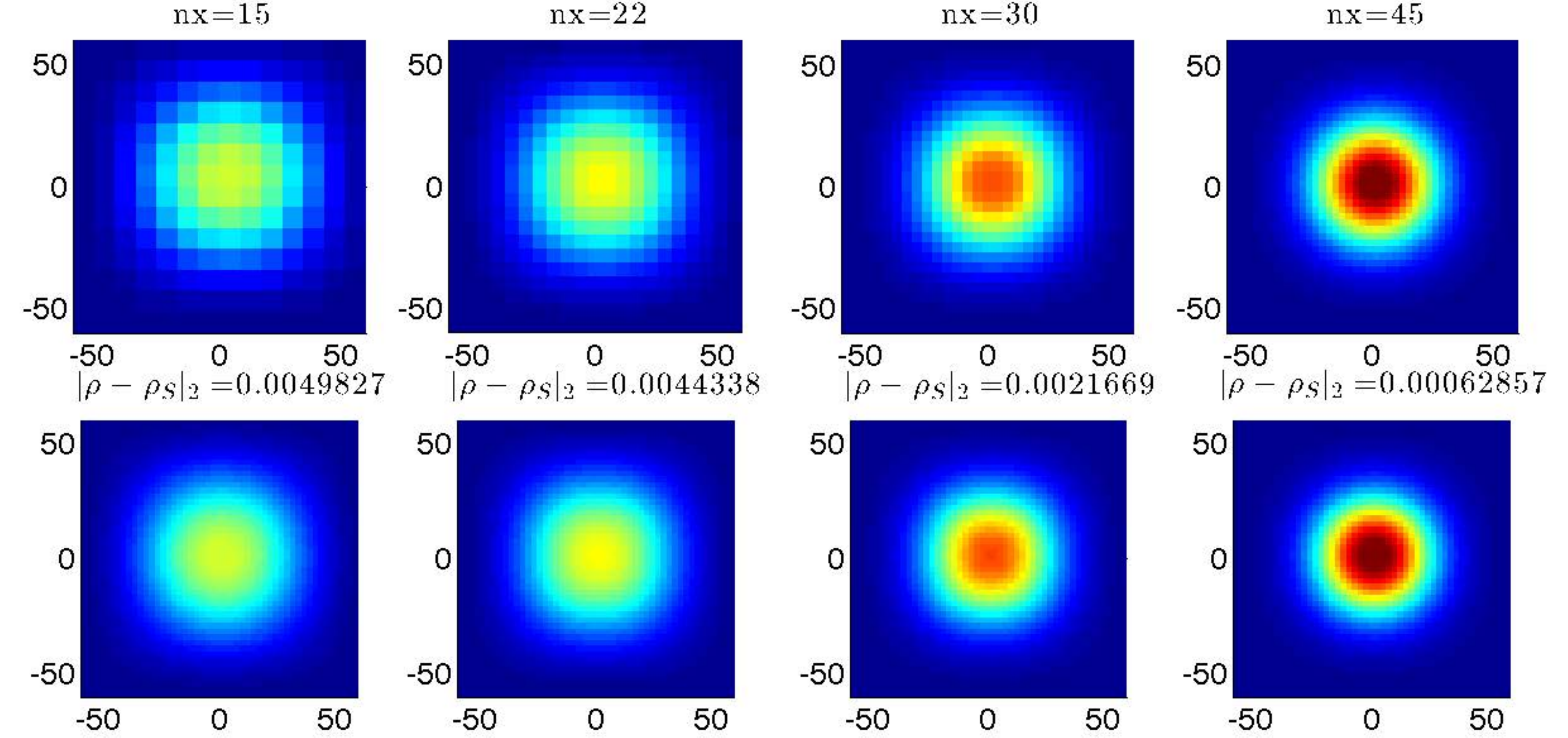}\\
\caption{The top row contains the original data of the numerical solution for $A=3.0$, below there is the interpolated solution, which is then compared to $\rho_R$. The colour scaling is the same as for the reference density in figure \ref{rhoDradial}. \label{rhoDdens}}
\end{figure}

This equation has been obtained in the previous section in the limit $A\rightarrow\infty$. For $A=3.0$ the stationary density is well approximated by the limit density. The computational domain is given by $x\in[-60,60]^2,v\in[-3,3]^2$. Table \ref{tab1} shows the number of grid points in each of the 4 dimensions, the grid size $h$ in the spatial domain, and the error.

\begin{figure}[ht]
\centering
\includegraphics[width=0.49\textwidth]{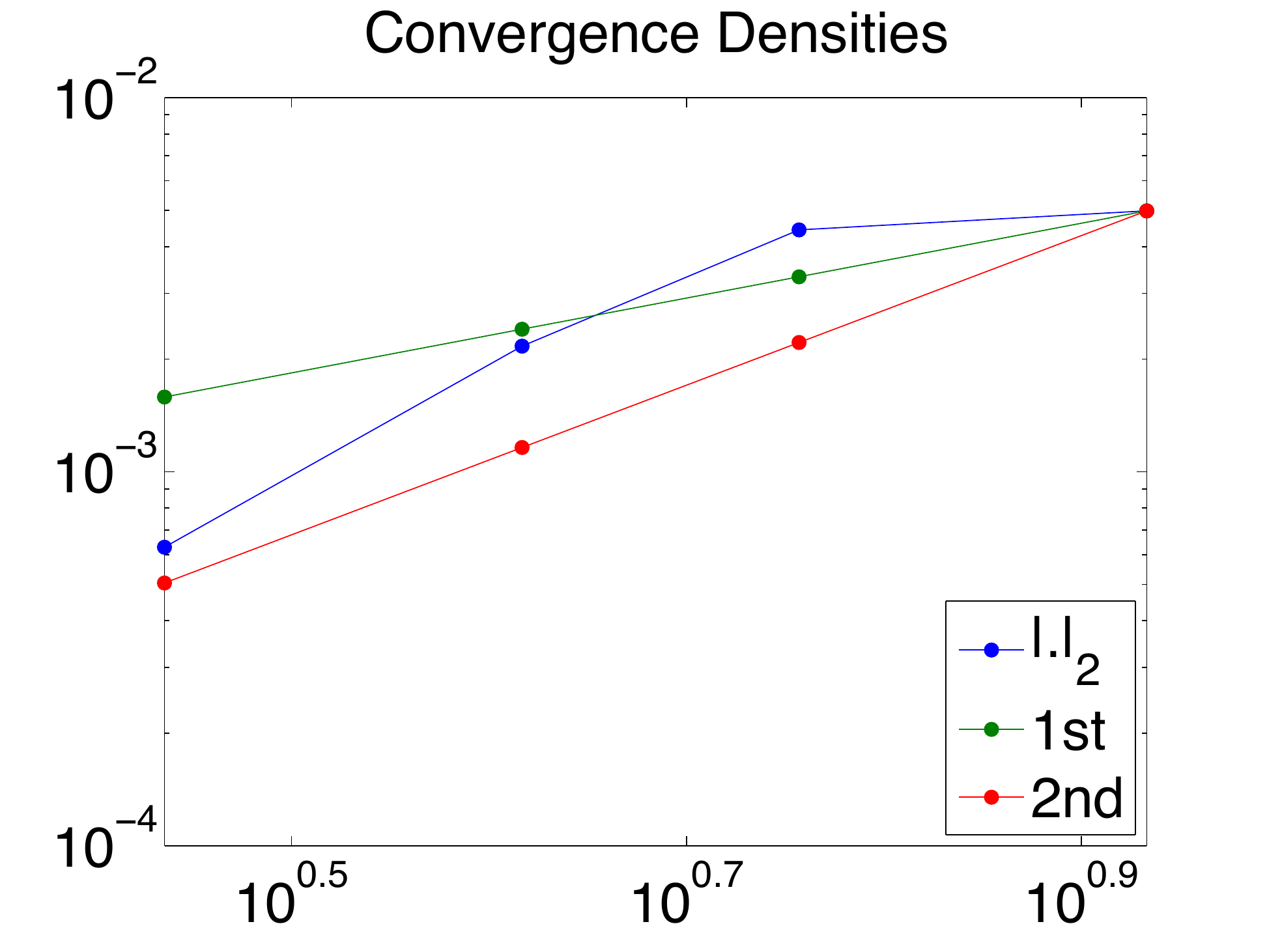}
\caption{The blue line shows the second order accuracy of the numerical scheme in $|\cdot|_2$ for $A=3.0$. The green and red lines are references for convergence rates of first or second order, respectively. \label{rhoDerror}}
\end{figure}

\begin{figure}[ht]
\centering
\includegraphics[width=0.47\textwidth]{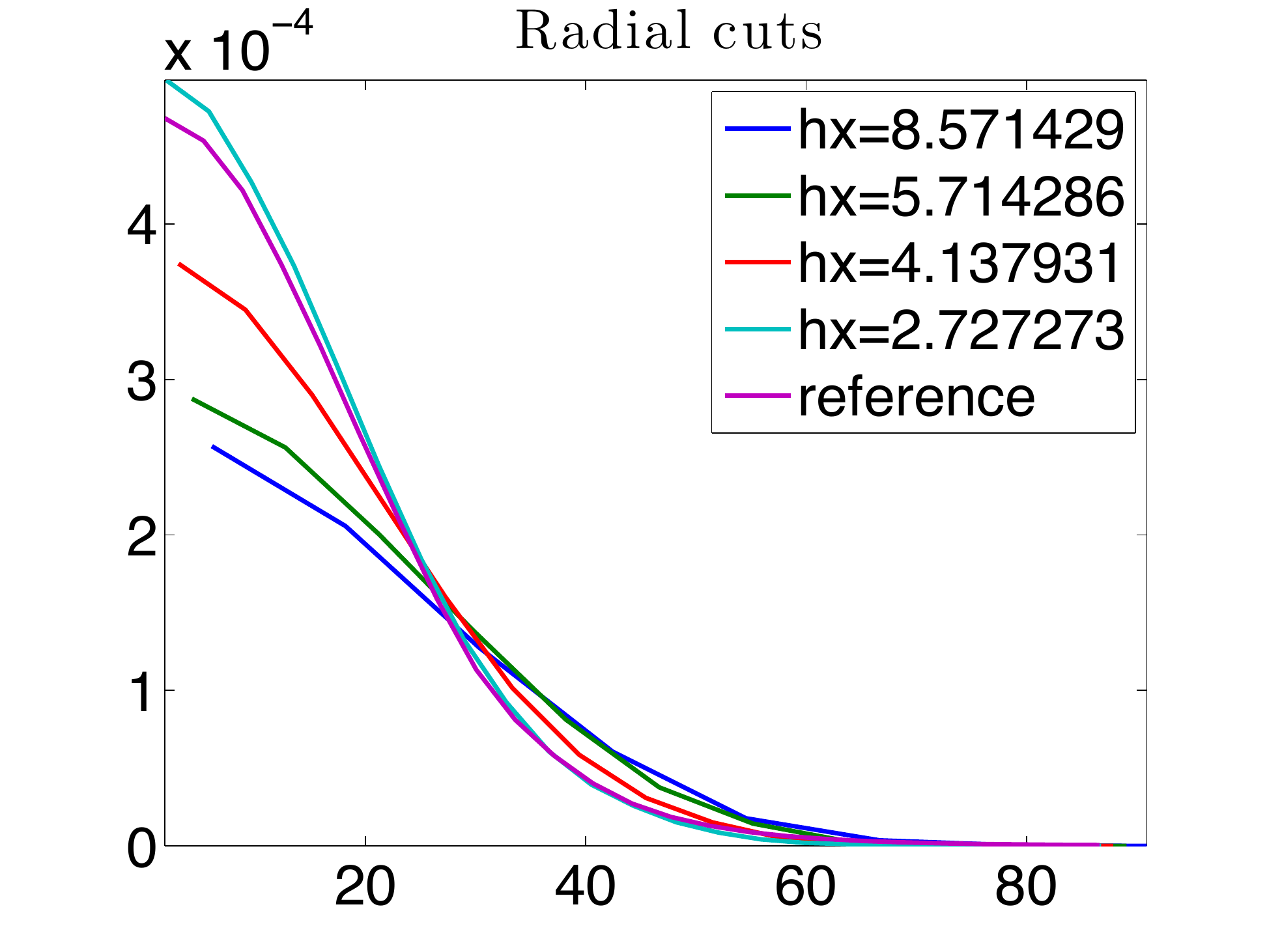}
\includegraphics[width=0.47\textwidth]{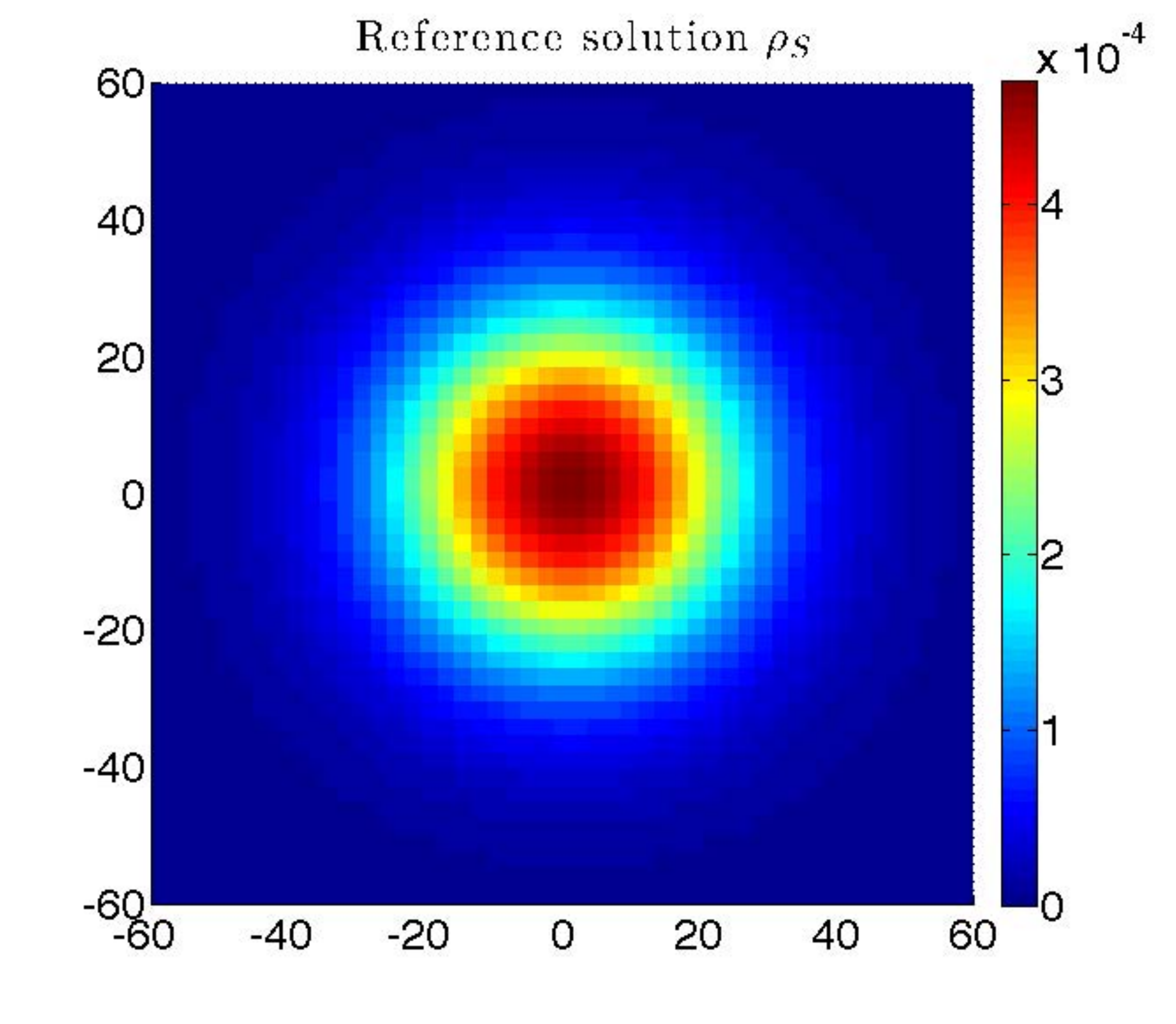}\\
\caption{On the left: radial plots of the numerical densities for the different grid sizes and the reference density $\rho_R$for $A=3.0$. On the right: the reference solution $\rho_R$ on its native grid.\label{rhoDradial}}
\end{figure}
\begin{table}[ht!]
\begin{center}
\begin{tabular}{|c|c|c|c|c|}\hline
$n$ & 15 &   22  &  30  &  45\\\hline
$h$ & 8.5714  &  5.7143 &   4.1379   & 2.7273\\\hline
$|\rho-\rho_R|_2$& 0.0049827 & 0.0044338 & 0.0021669 & 0.0008489 \\ \hline
\end{tabular}
\end{center}
\caption{Grid size and error for $A=3.0$.}
\label{tab1}
\end{table}

The  step sizes in the velocity domain are  decreased in the same way as the spatial grid $h$ for the convergence plots. In fact, we use as many grid points for the spatial as for the velocity grid. The respective densities can be found in figure \ref{rhoDdens}. We get the logarithmic error plots in figure \ref{rhoDerror}. Radial plots for the densities and the reference density can be found in figure \ref{rhoDradial}. Note, that the colour scaling is the same for reference density in figure \ref{rhoDradial} and the densities in figure \ref{rhoDdens}.

\subsubsection{A=0.0}
In this case without diffusion, we  obtain a single mill solution as equilibrium solution. Since there is no  analytical solution to compare,  we take as a reference solution the numerical solution with a
fine resolution. In this investigation, we use the solution with 60 grid points in each of the 4 directions as the reference solution.  The computational domain will be given by $x\in[-50,50]^2,v\in[-3,3]^2$. We  denote the density of the fine resolved numerical solution by $\rho_R$ and  compare them to the densities $\rho$ obtained from the  lower resolutions, in the norm stated in (\ref{norms}).
Table \ref{tab2} shows the number of grid points in each of the 4 dimensions, the grid size $h$ in the spatial domain, and the error.

\begin{figure}
\centering
\includegraphics[width=\textwidth]{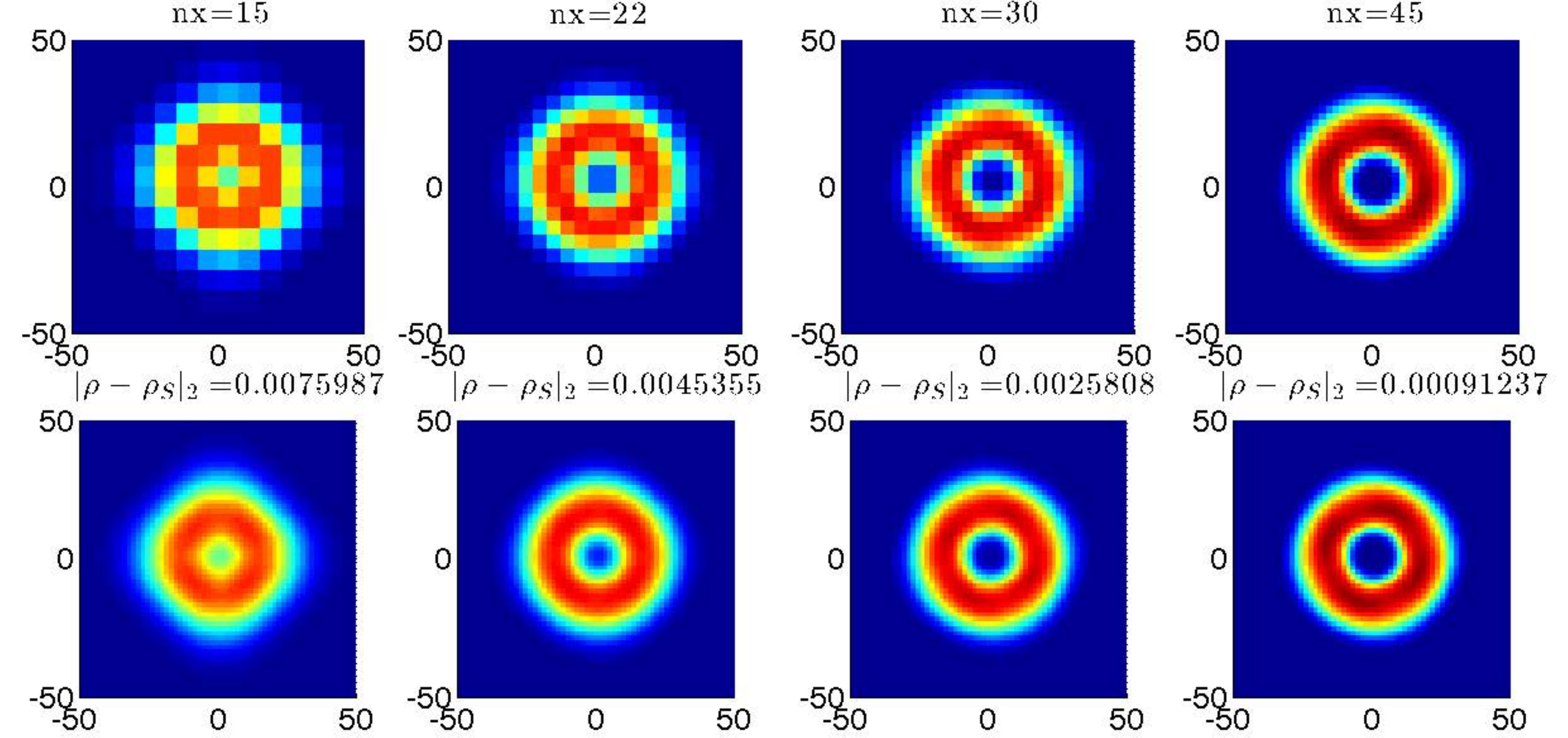}
\caption{The top row contains the original data of the numerical solution for $A=0.0$, below there is the interpolated solution, which is then compared to $\rho_R$. The color scaling is the same as for the reference density in figure \ref{rhoSradial}.\label{rhoSdens}}
\end{figure}

\begin{figure}
\centering
\includegraphics[width=0.49\textwidth]{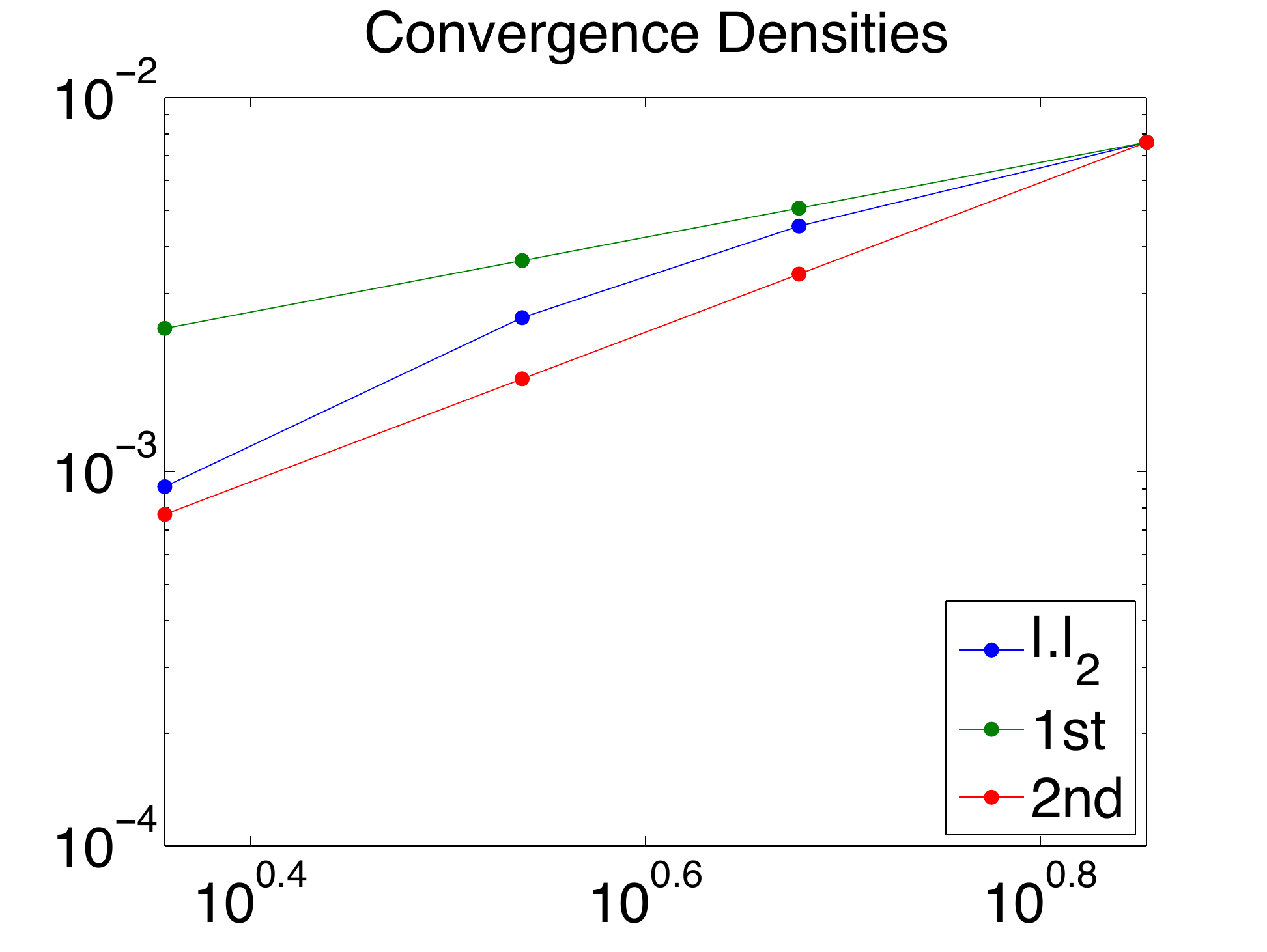}
\caption{The blue line shows the errors of the numerical scheme in $|\cdot|_2$ for $A=0.0$. The green and red lines show the errors corresponding to a convergence rate of first or second order, respectively. Since the error line is parallel to the red line, we conclude that the numerical scheme is of second order for $A=0.0$.\label{rhoSerror}}
\end{figure}

\begin{table}[ht!]
\begin{center}
\begin{tabular}{|c|c|c|c|c|}\hline
$n$ & 15  &  22 &   30 &   45\\\hline
$h$ & 7.1429  &  4.7619  &  3.4483 &   2.2727\\\hline
$|\rho-\rho_R|_2$ & 0.0075987 & 0.0045355 & 0.0025808 & 0.0009124 \\\hline
\end{tabular}
\end{center}
\caption{Grid size and error for $A=0$.}
\label{tab2}
\end{table}

\begin{figure}
\centering
\includegraphics[width=0.47\textwidth]{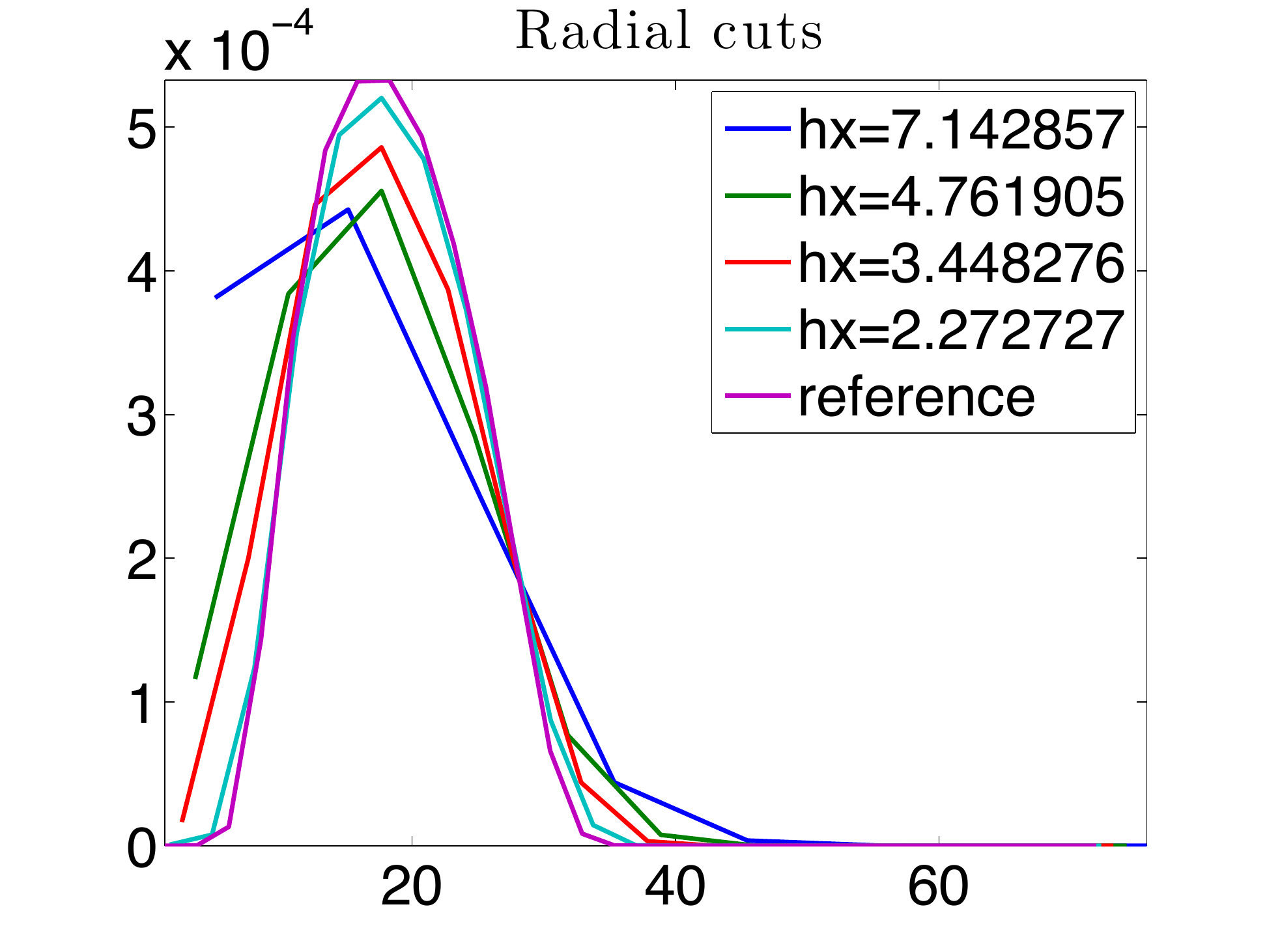}
\includegraphics[width=0.47\textwidth]{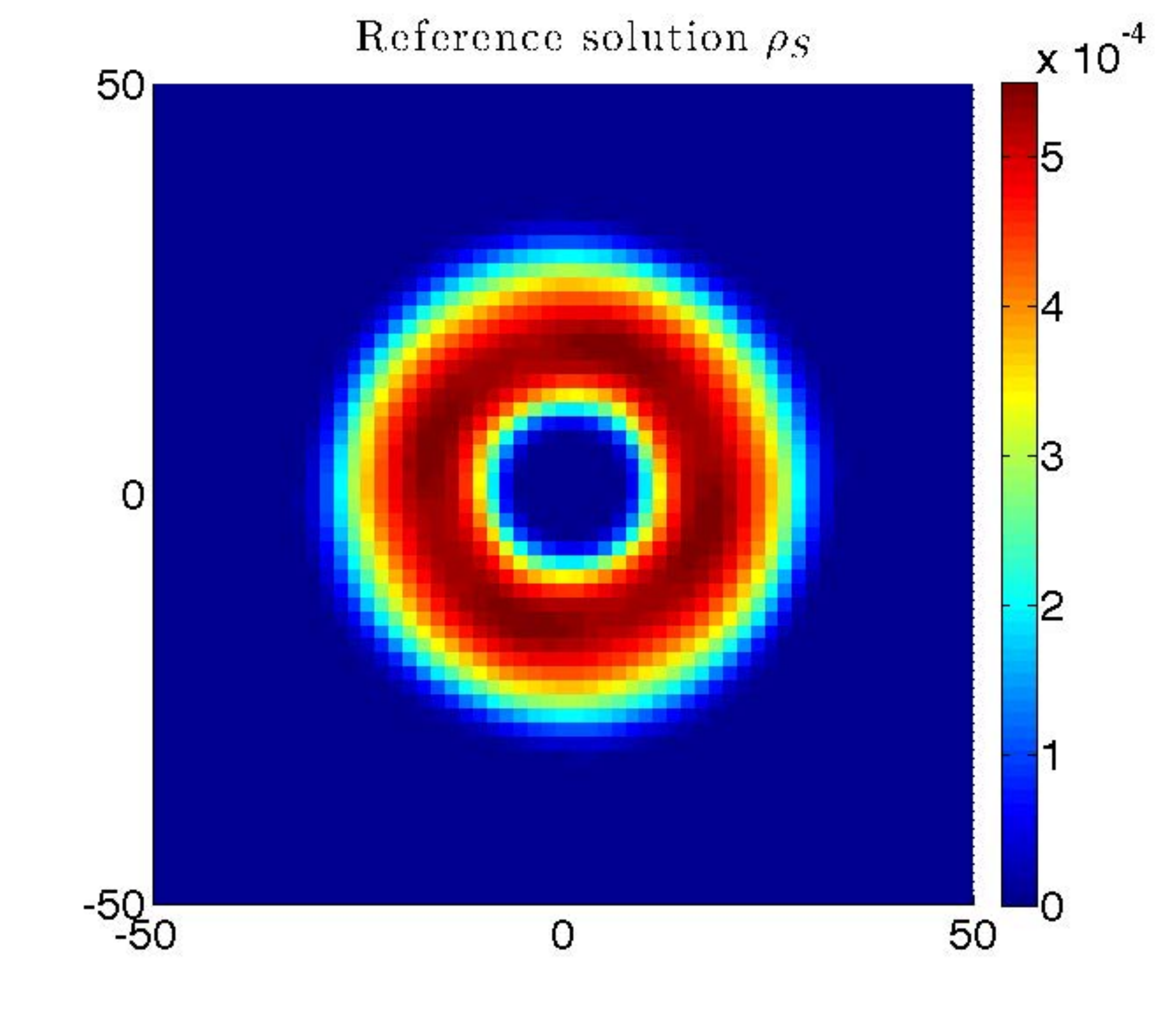}
\caption{On the left: radial plots of the numerical densities for the different grid sizes and the reference density $\rho_R$ for $A=0.0$.  On the right: the reference solution $\rho_R$ on its native grid.\label{rhoSradial}}
\end{figure}

The respective densities can be found in figure \ref{rhoSdens}. From this, we get the logarithmic error plot in figure \ref{rhoSerror}. Radial plots for the densities and the reference solution can be found in figure \ref{rhoSradial}.

\section{Numerical investigation of milling States} 
\label{Numerics}
In this section we look at numerical solutions for the full problem (\ref{fpstandard}) and compare them to the microsopic solutions of equation (\ref{micro}) and the macroscopic approximations derived in the last section. For $\alpha, \beta$ we use the same numerical values as in the last section. Moreover, the constants of the interaction potential are chosen as in \cite{hauptpaper} as  $C_a=20,C_r=50,l_a=100,l_r=2$. Depending on the initial conditions  special stationary solutions, so called   single or  double mills, appear for diffusion coefficient $A=0$. They are obtained as  stationary states of the interacting particle system as well as   the kinetic equation. We discuss these stationary states numerically and concentrate on the investigation of  the behaviour of these solutions for increasing noise parameter $A$.

\subsection{Milling solutions}
\begin{figure}
    \centering
    \includegraphics[width=6.3cm]{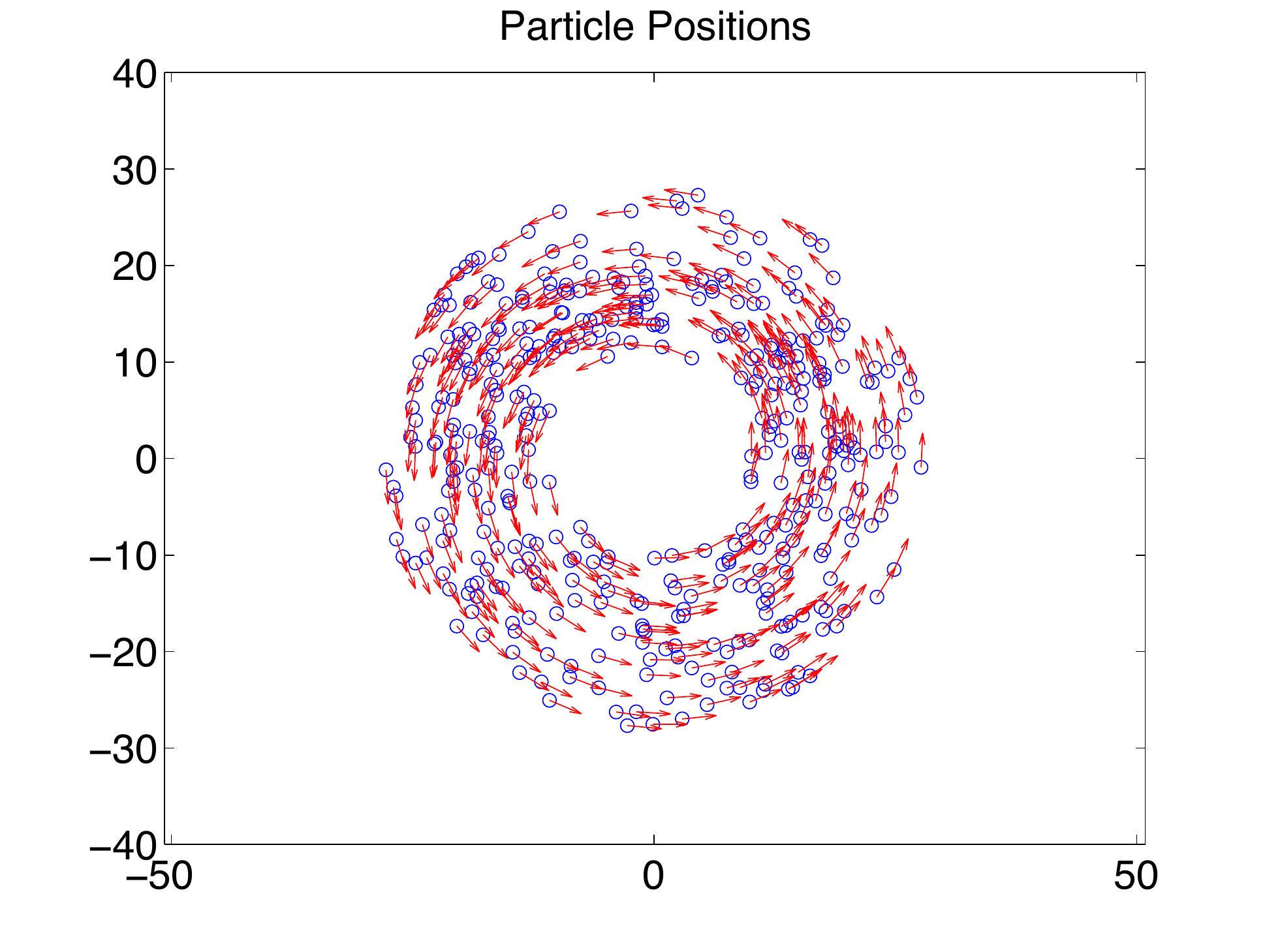}
    \includegraphics[width=6.3cm]{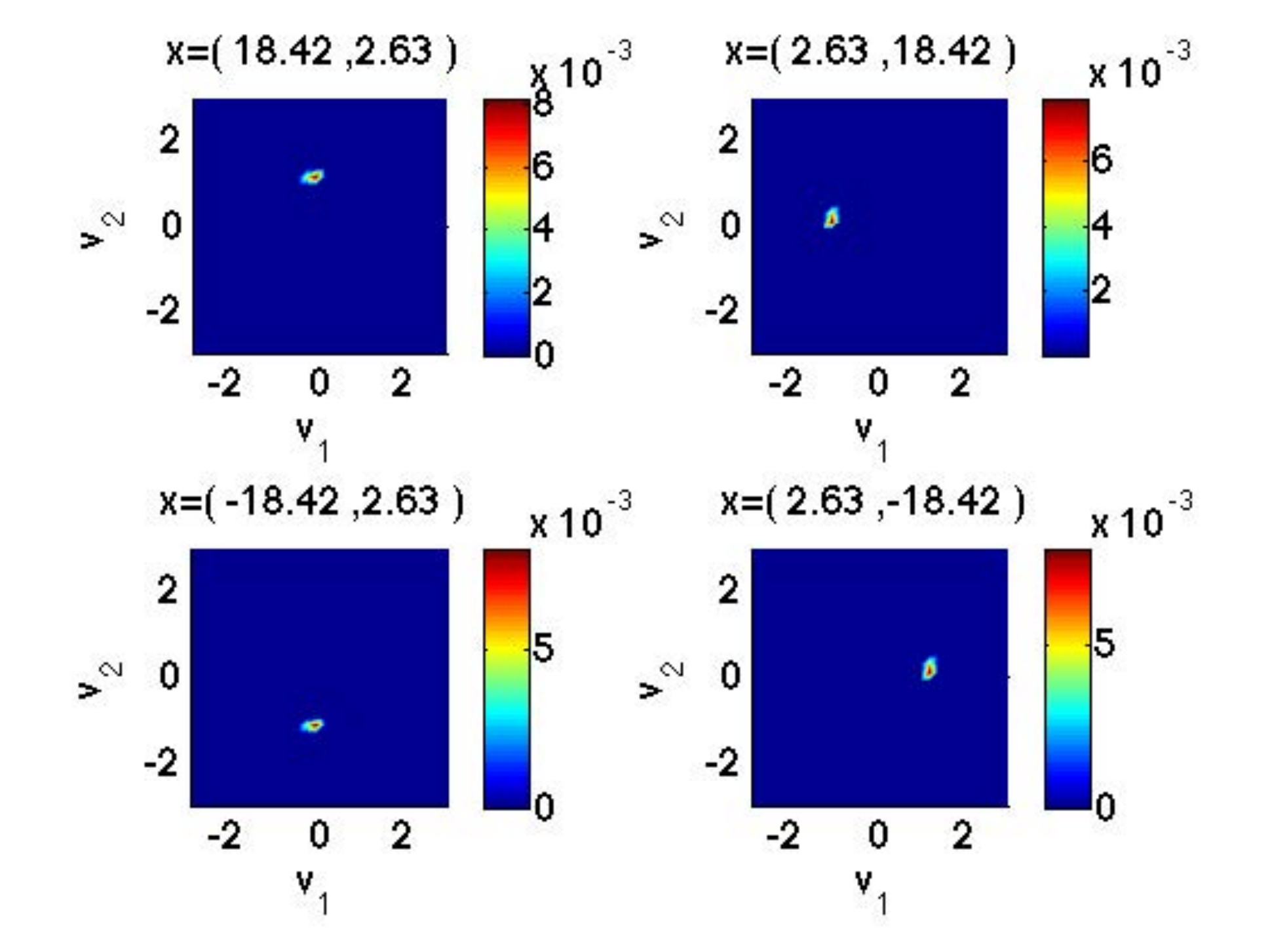}\\
    \includegraphics[width=6.3cm]{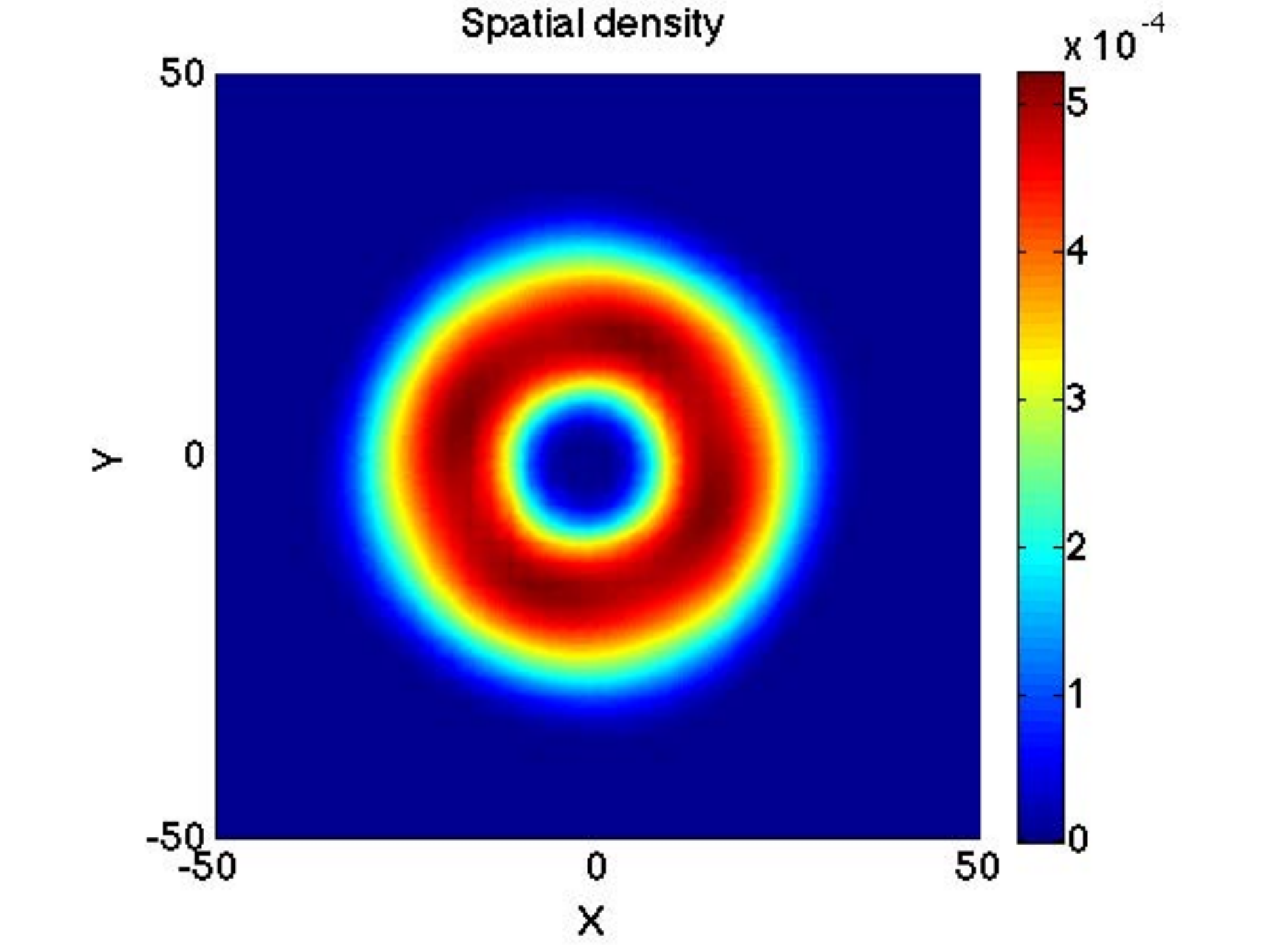}
    \includegraphics[width=6.3cm]{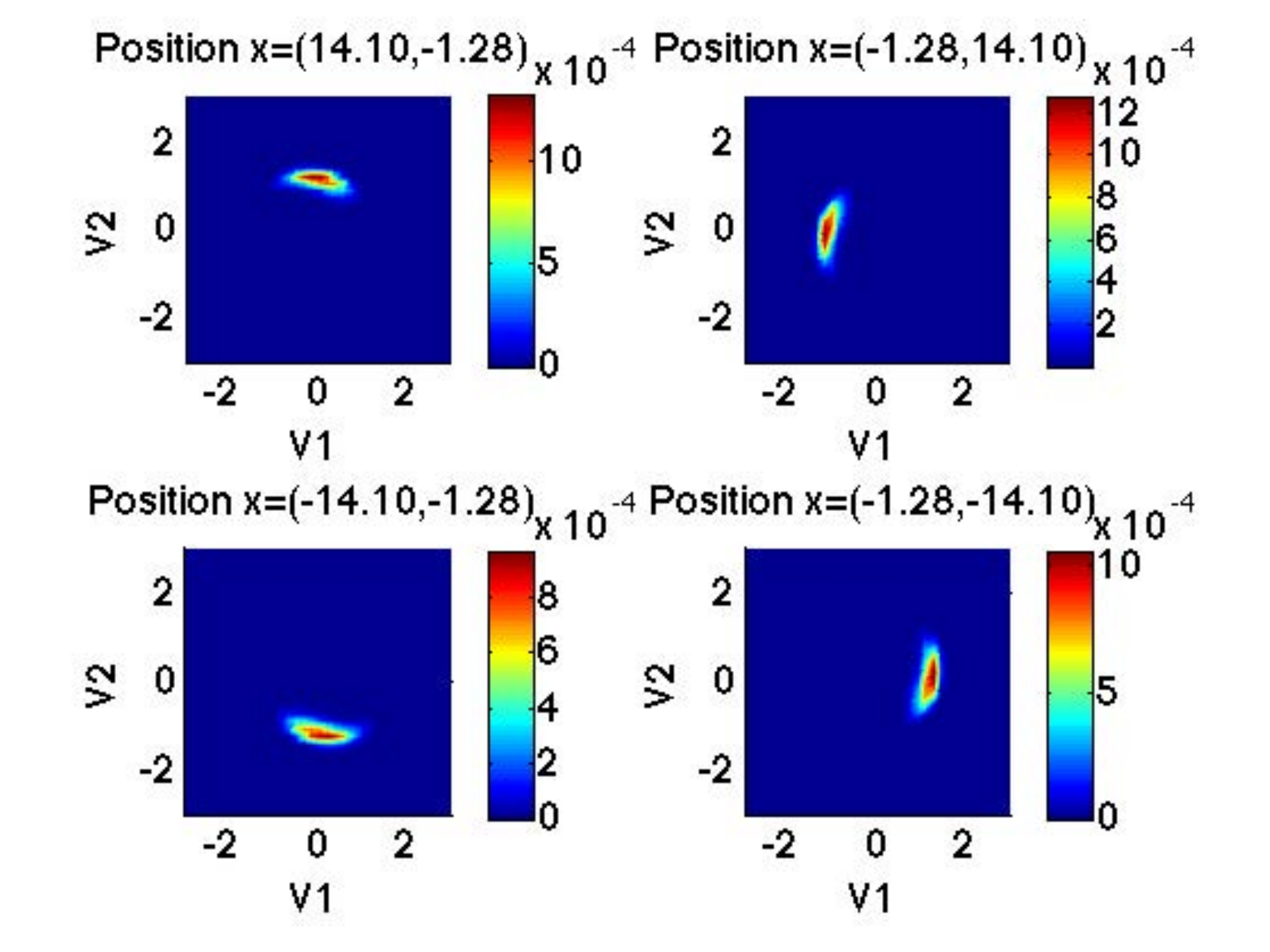}
    \caption{Comparison of single mills from microscopic equations (top) and from kinetic equations (bottom). To the left, we have the particle positions and kinetic particle distribution, to the right, there are the velocity distributions at fixed points on the grid. Note, that the grid in velocity space is different for the microscopic and kinetic case, so the positions do not coincide exactly.}\label{singlemillkinmiccomp}
  \end{figure}

The single mill solution in the microscopic context is a structure, where all the particles with positions $x_i\in\mathbb{R}^2$ go around the origin with velocities 
$$
v_i=\sqrt{\frac{\alpha}{\beta}}\frac{x_i^\perp}{||x_i||}.
$$ 
In the kinetic case we display  the spatial distribution of the density $\rho=\int f\;dv$. The velocities are visualized by looking at the velocity distributions $f(x_0,v)$ at fixed spatial points $x_0\in\operatorname{supp}\rho$. If a single mill is obtained,  $f(x_0,v)$ 
will become concentrated at 
$$
v\sim \sqrt{\frac{\alpha}{\beta}}\frac{x_0^\perp}{||x_0||}.
$$ 
For comparison, we  generate a distribution function $f$ from the microscopic results,  by counting particles in a 4D-histogram, which is then normalized to mass 1 as for the kinetic solution. We note, that a large amount of particles is needed to generate a reasonably smooth  histogram. In figure \ref{singlemillkinmiccomp}, we plotted 400 microscopic particle positions with arrows indicating their velocities. The actual computation was carried out with 13000 interacting particles, where we have chosen 400 particles randomly for plotting the results. The histogram is based on the position and velocity data of all 13000 particles. The spatial density obtained from the kinetic equation and pictures of the velocity distributions at corresponding fixed points obtained from microscopic and kinetic solution are also shown in figure \ref{singlemillkinmiccomp}. The results coincide  well, although the peaks in the kinetic velocity distribution are broader than the microscopic ones, due to the diffusivity of the kinetic solver.

The single mill is obtained for $A=0$ for a specific type of initial conditions. One has to prepare the system in a state which is not very far from the single mill solution. In the kinetic case, this can be realised by using the initial condition
\begin{align*}
 f_0(x,v)&=\begin{cases}
 C&\mbox{, if }12\leq||x||\leq 29, \; \sqrt{\frac{\alpha}{\beta}}-\frac{1}{2}\leq||v||\leq\sqrt{\frac{\alpha}{\beta}}+\frac{1}{2}, \; ||\frac{v}{||v||}-\frac{x^\perp}{||x||}||\leq \frac{15}{100}.\\
 0&\mbox{, else}.
 \end{cases}
 \end{align*}
$C$ is chosen in such a way, that $\int f_0\;dxdv=1$. One can achieve a corresponding initial condition in the microscopic case by placing particles $i$ randomly at positions $x_i$ with $12\leq||x_i||\leq 29$ and assigning random velocities $v_i$ with the corresponding deviation from  the direction $x_i^\perp$.

 \begin{figure}[ht!]
  \centering
  \includegraphics[width=6.3cm]{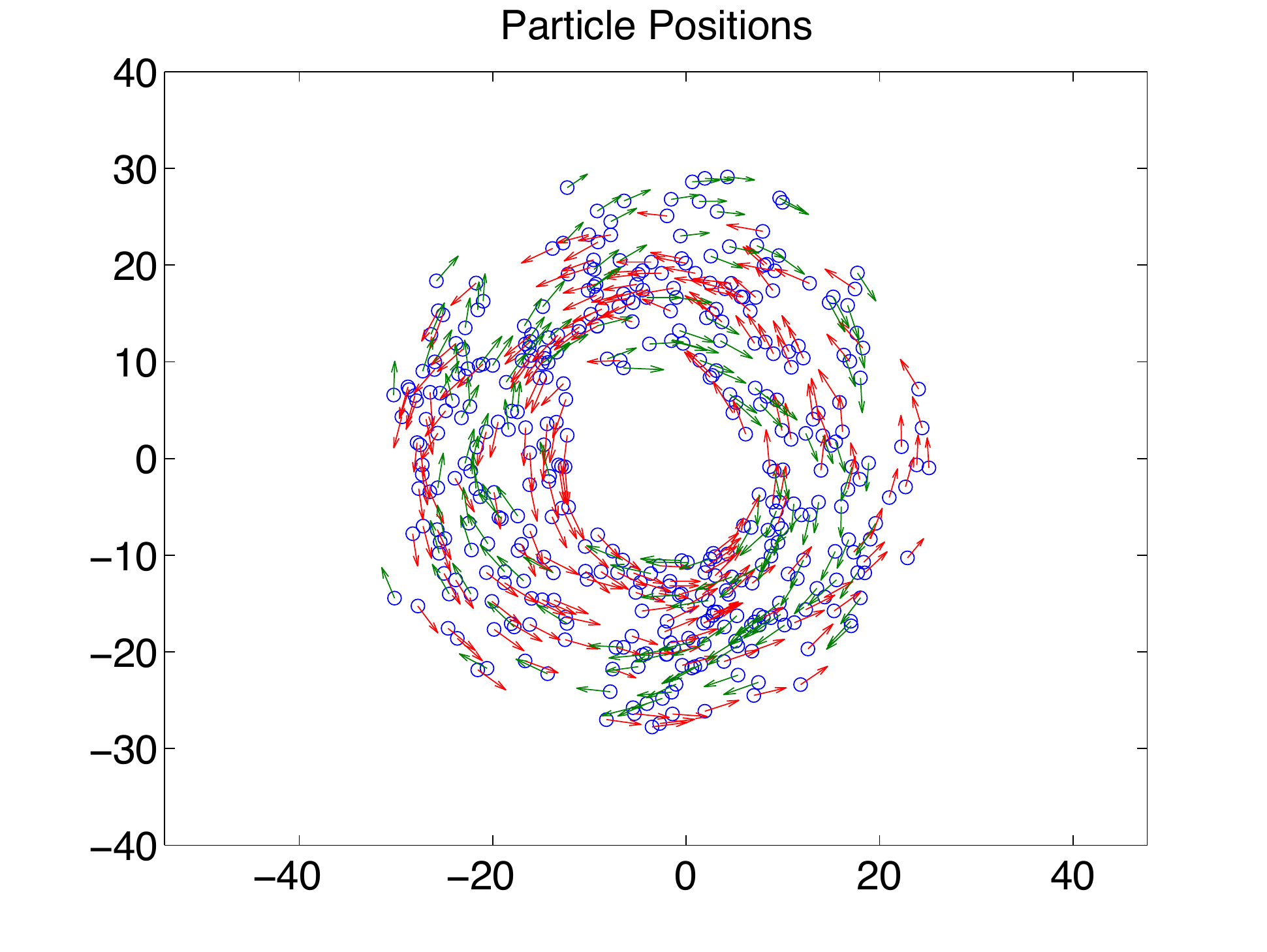}
  \includegraphics[width=6.3cm]{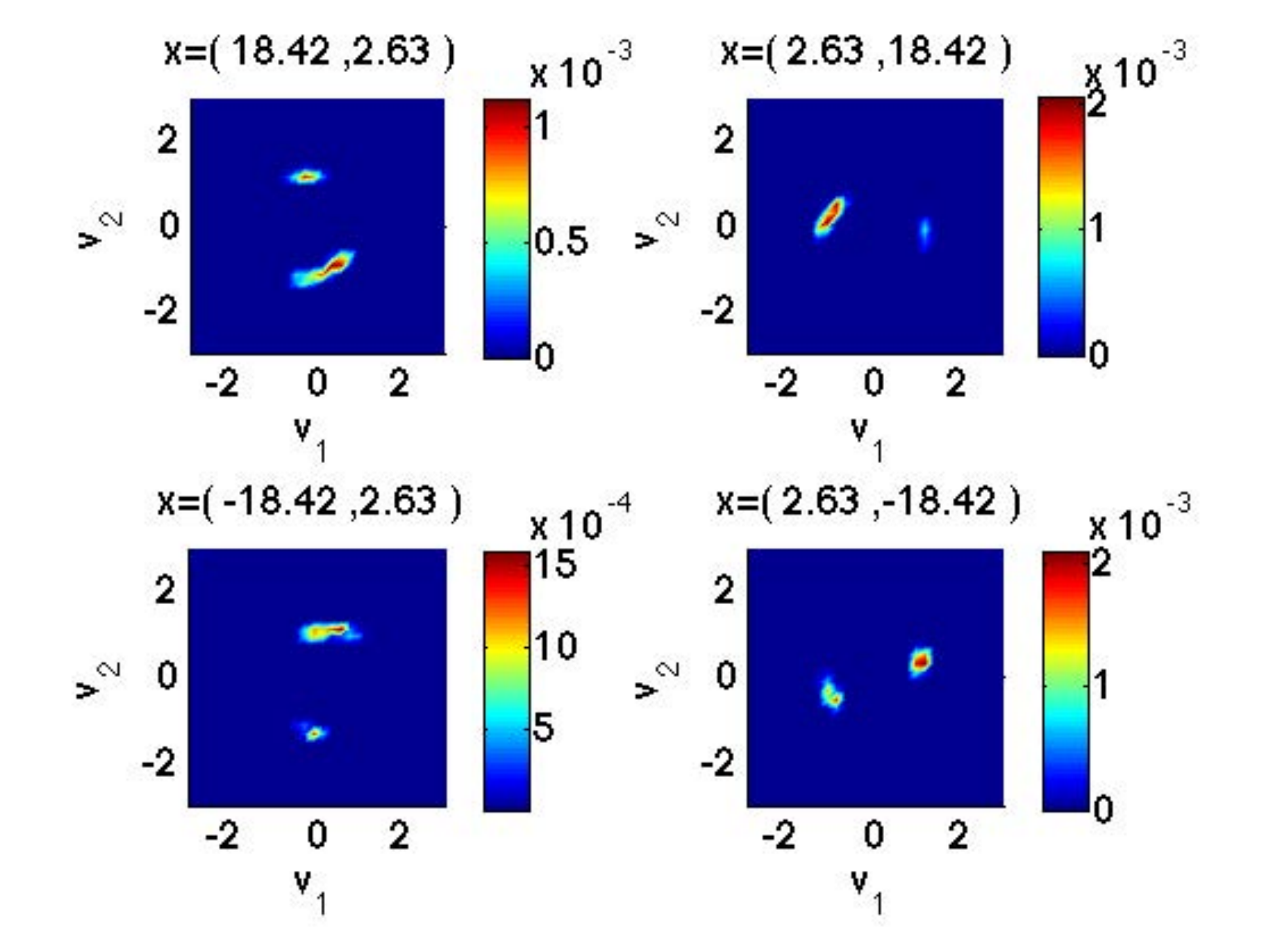}\\
  \includegraphics[width=6.3cm]{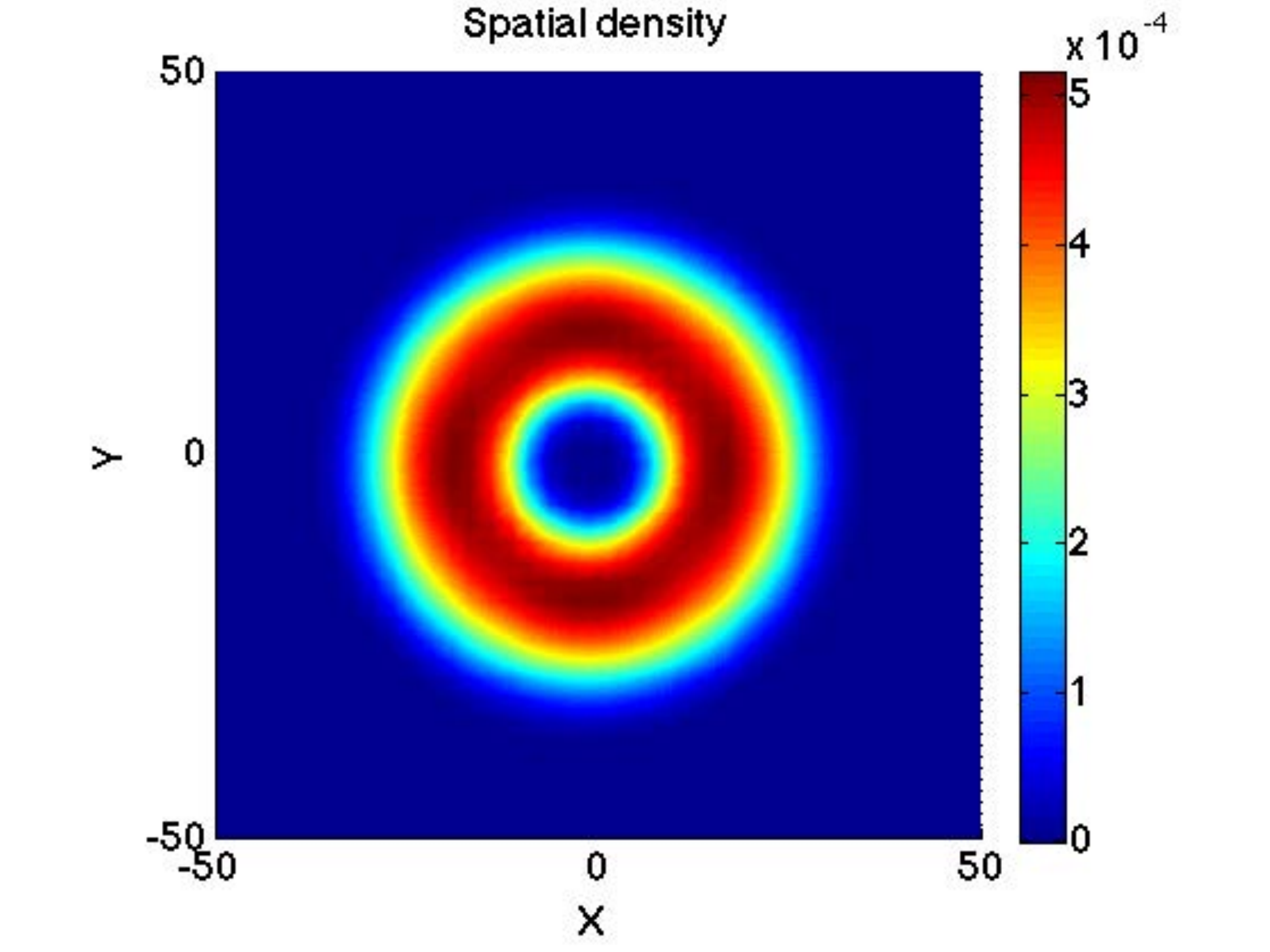}
  \includegraphics[width=6.3cm]{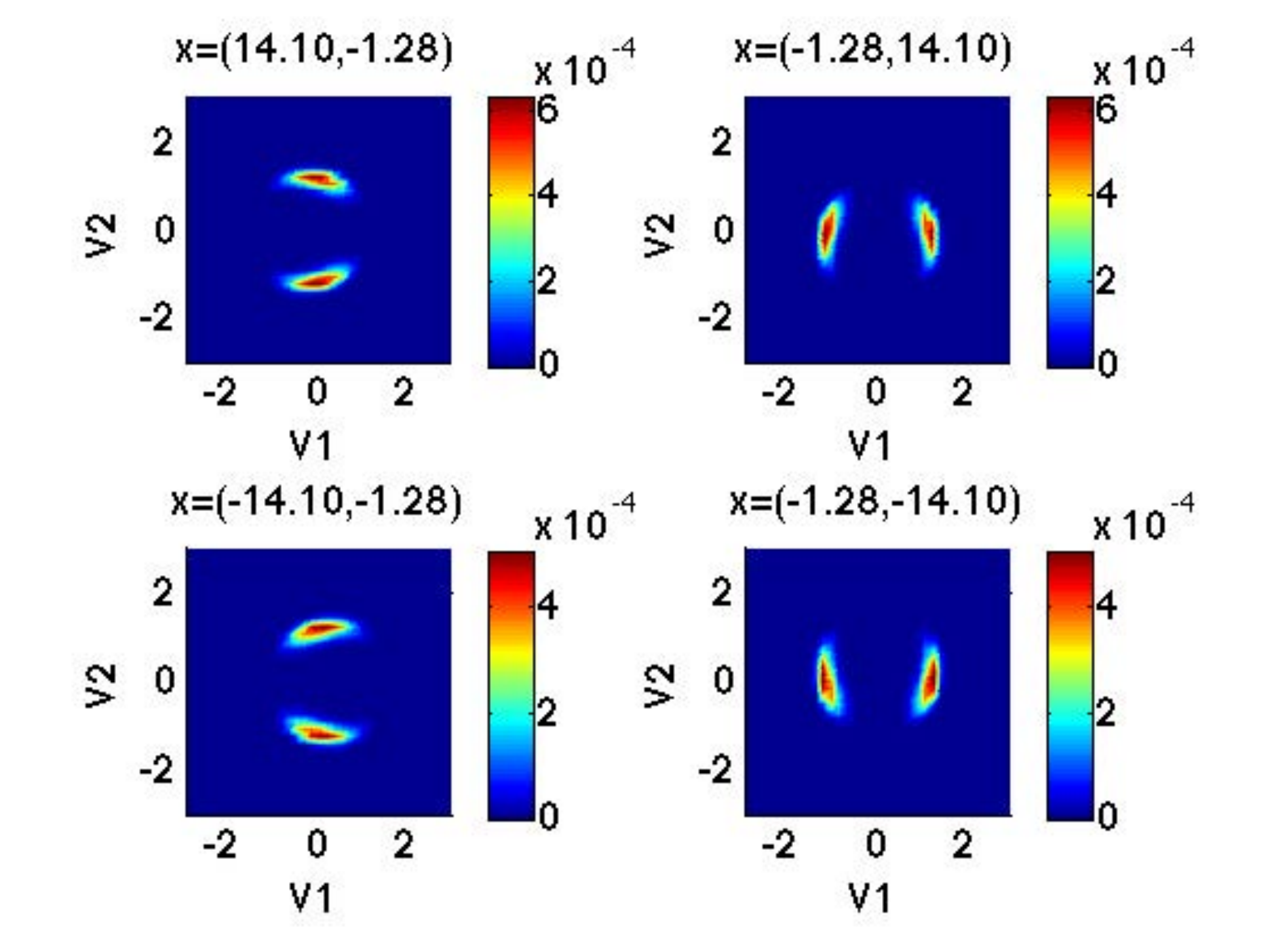}
  \caption{Comparison of double mills from microscopic equations (top) and from kinetic equations (bottom). To the left, we have the particle positions and kinetic particle distribution, to the right, there are the velocity distributions at fixed points on the grid. Note, that the grid in velocity space is different for the microscopic and kinetic case, so the positions do not coincide exactly.}\label{doublemillkinmiccomp}
  \end{figure}

A double mill is, from the  microscopic point of view, a situation, where one part of the particles at the 
positions $x_i$ is circulating around the origin with speed 
$$
v_i=\sqrt{\frac{\alpha}{\beta}}\frac{x_i^\perp}{||x_i||},
$$
while the other part at around the same location $x_i$ is going in the opposite direction 
$-v_i$. The kinetic spatial density looks like the one for the single mill solution, but the velocity distributions at fixed points show two symmetric peaks in $f$ at the sphere $|v|^2=\tfrac{\alpha}{\beta}$. Computing a histogram from the microscopic solution as before, we compare the results in figure \ref{doublemillkinmiccomp}. Again, the two solutions agree reasonably well. Note, that it is not possible to obtain a double mill with the macroscopic model because of the monokinetic closure, see \cite{CDP}.
  
The double mill is again obtained for $A=0$ by choosing appropriate initial conditions. For example,
one may use the following distribution function for the kinetic equation:
\begin{align*}
 f_0(x,v)&=\begin{cases}
 C&\mbox{, if }12\leq||x||\leq 29, \; ||v||\leq 1.6\\
 0&\mbox{, else}.
 \end{cases}
 \end{align*}
 $C$ is chosen as before. The corresponding initial values  in the microscopic context is obtained by  placing particles randomly in the same region in space and assign random velocity vectors $v_i$, which are limited by $||v_i||\leq 1.6$. Note, that the initial condition to get a double mill pattern is less restrictive than the one for the single mill. This indicates that  a double mill is in a certain sense more stable or robust than a single mill for small stochasticity. This fact will become clearer investigating situations with increasing  noise parameter $A$ in the next subsection.

\begin{figure}[Hbt]
\centering
\includegraphics[width=6.2cm]{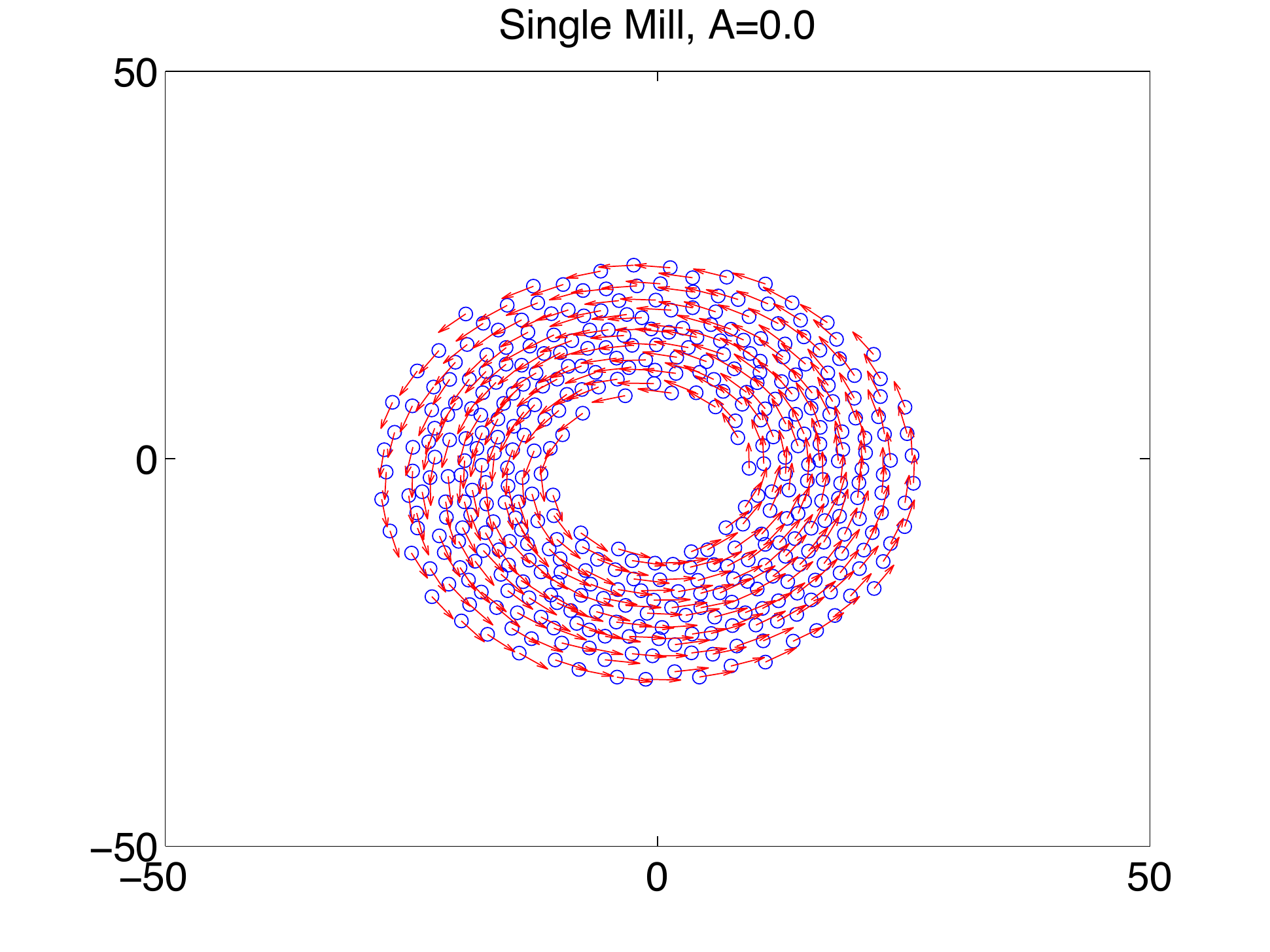}
\includegraphics[width=6.2cm]{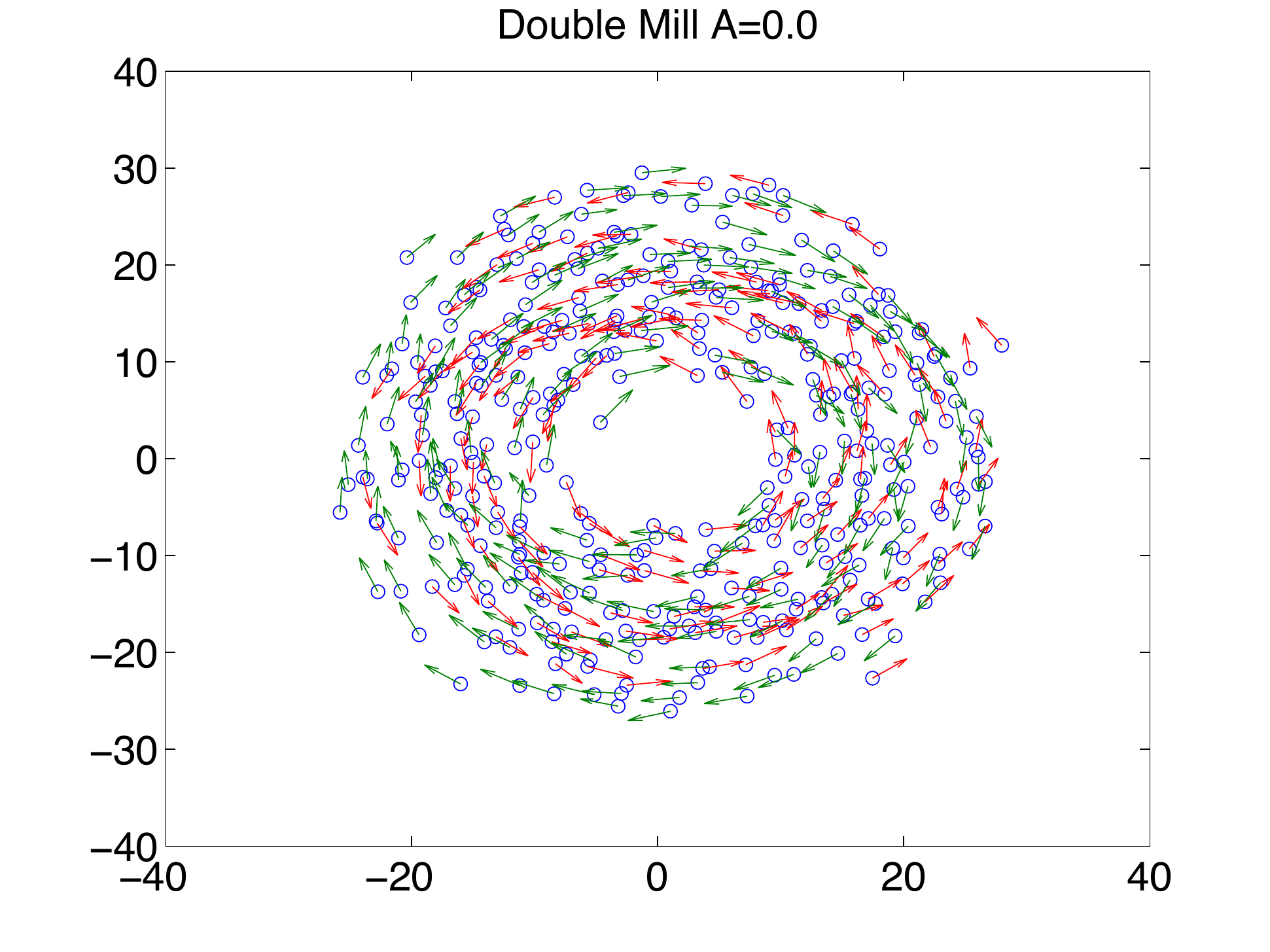}\\
\includegraphics[width=6.2cm]{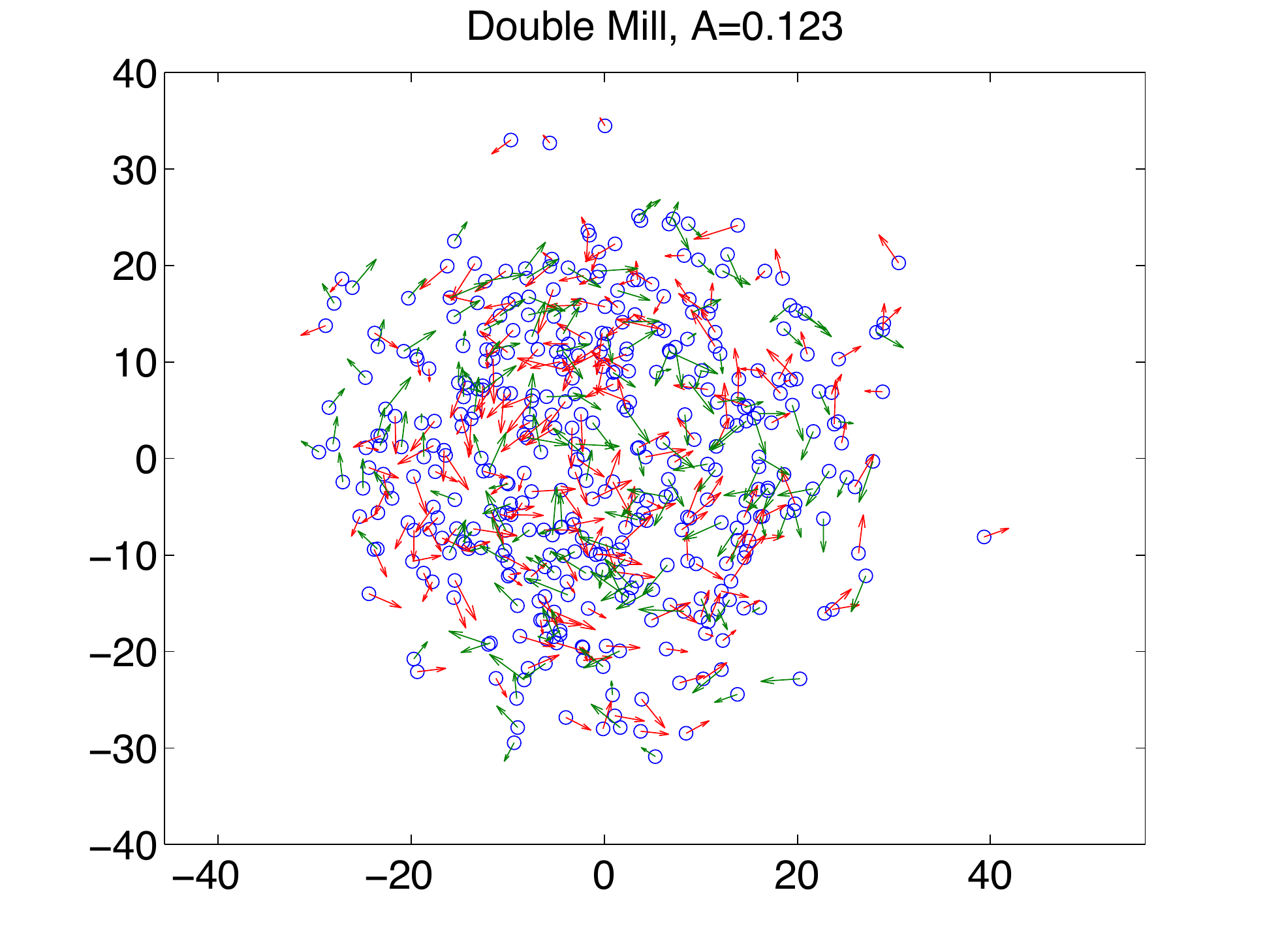}
\includegraphics[width=6.2cm]{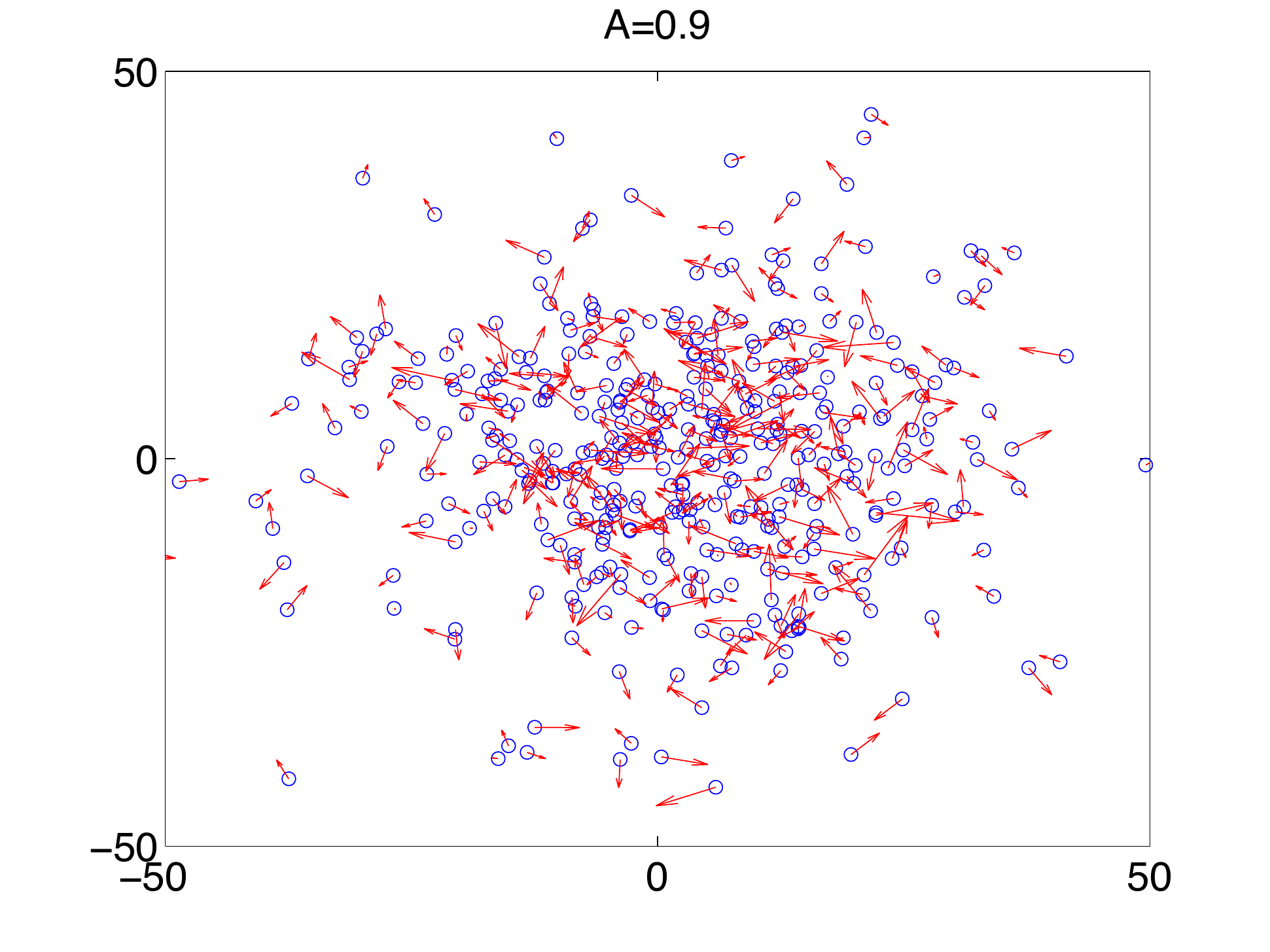}
\caption{Solutions of the microscopic system with 400 particles: At the top on the left, one can see a single mill for $A=0$, on the right, we have a double mill for $A=0$. Below, $A$ is raised first to 0.123, where we still have a double mill  until we end up in an unordered state at $A=0.9$.\label{millsmikro}}
\end{figure}

\subsection{Influence of the noise amplitude: Phase Transition}
In this section the microscopic and kinetic equations are investigated for  different values of $A$ starting  with a single and a double mill initial condition. Figure \ref{millsmikro} shows the behaviour of the microscopic system for three values of the noise coefficient: For $A=0$, we obtain the expected single mill or double mill depending on the initial conditions, then for $A=0.123$ the stationary state is still a double mill independent of the initial state being a single or a double mill initial distribution. For $A=0.9$, we have an unordered state.

\begin{figure}
\centering
\includegraphics[width=6.3cm]{kineticSpatialDensA0000.pdf}
\includegraphics[width=6.3cm]{kineticVelocitiesA0000.pdf}\\
\includegraphics[width=6.3cm]{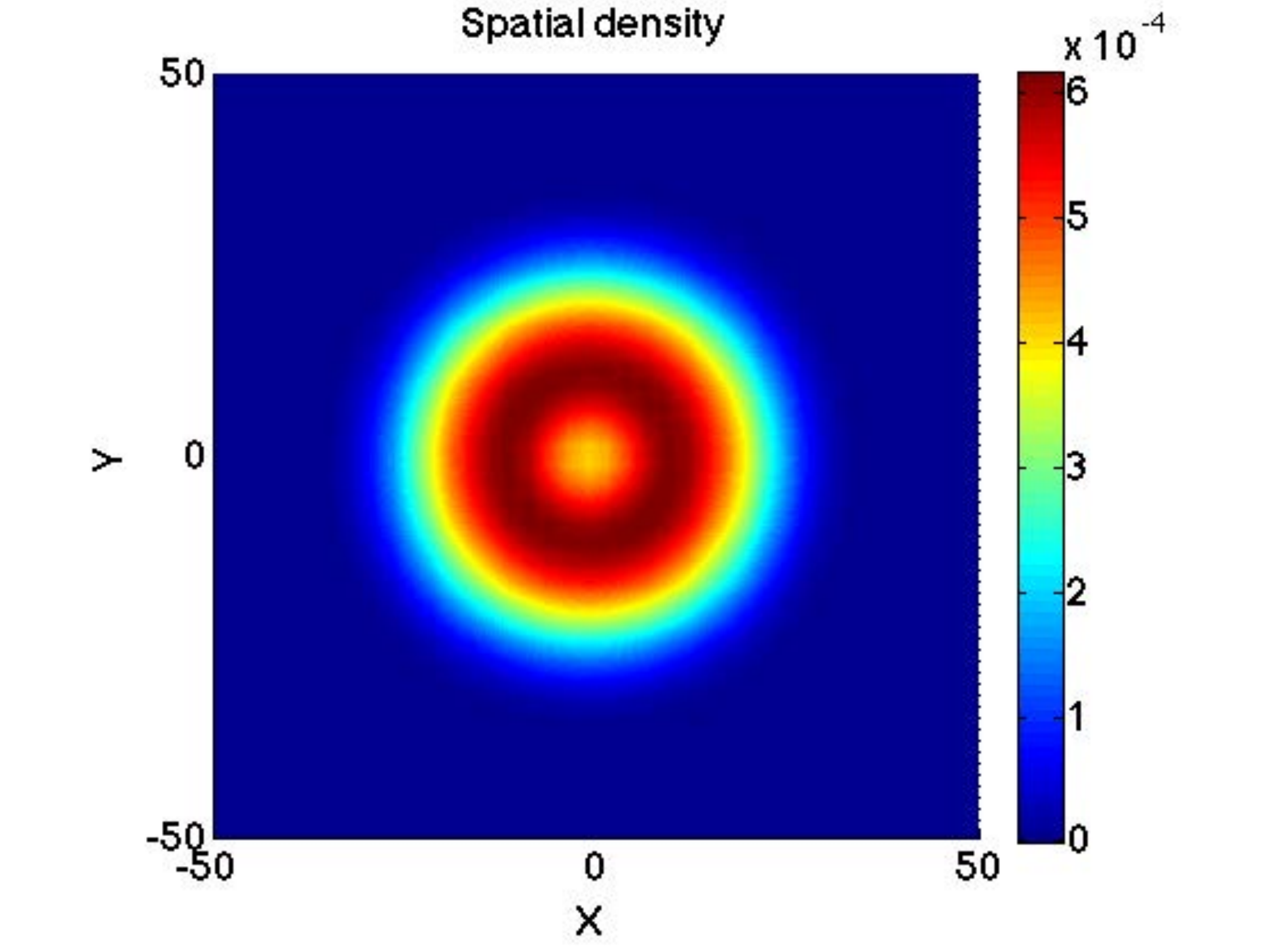}
\includegraphics[width=6.3cm]{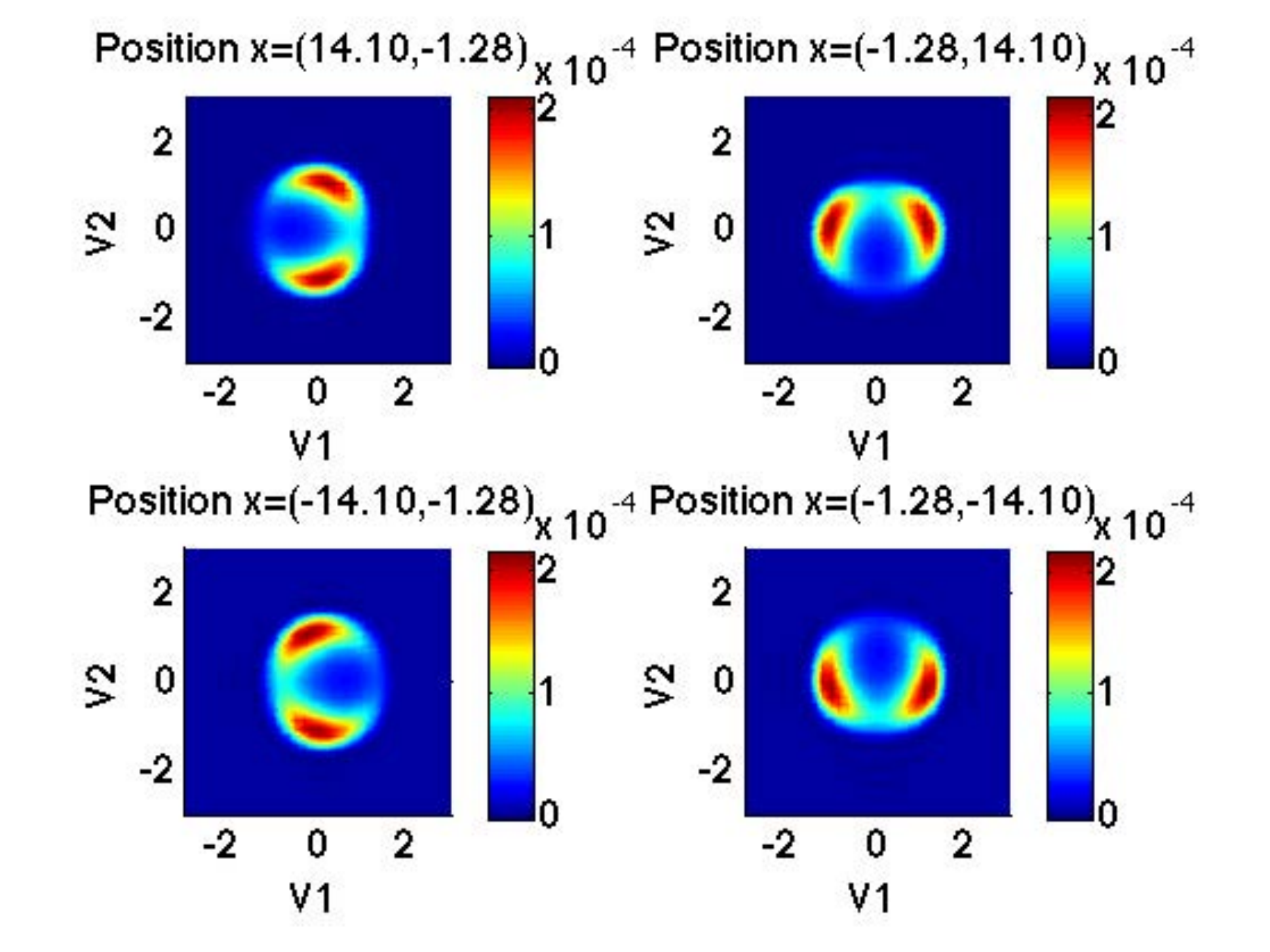}\\
\includegraphics[width=6.3cm]{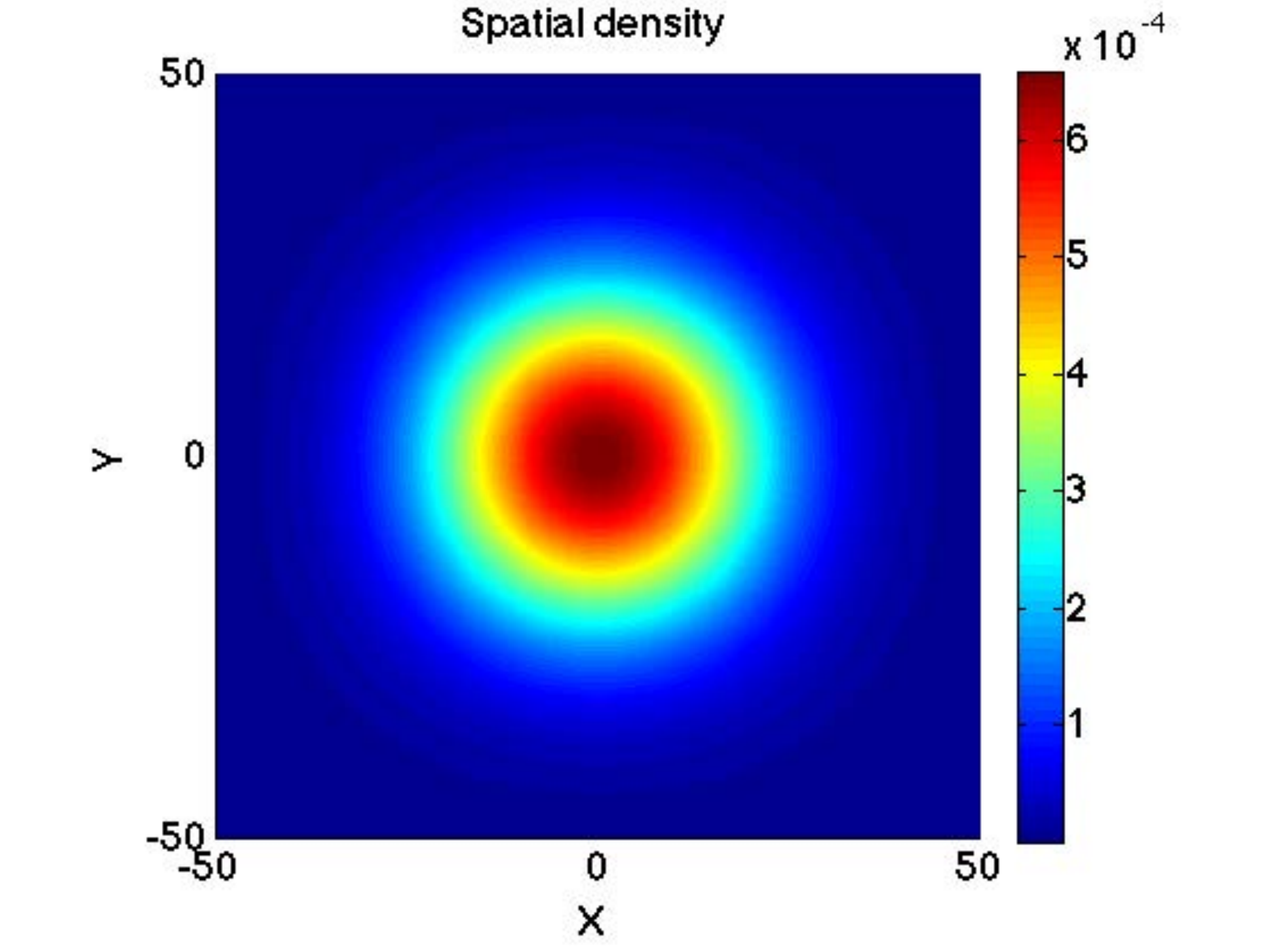}
\includegraphics[width=6.3cm]{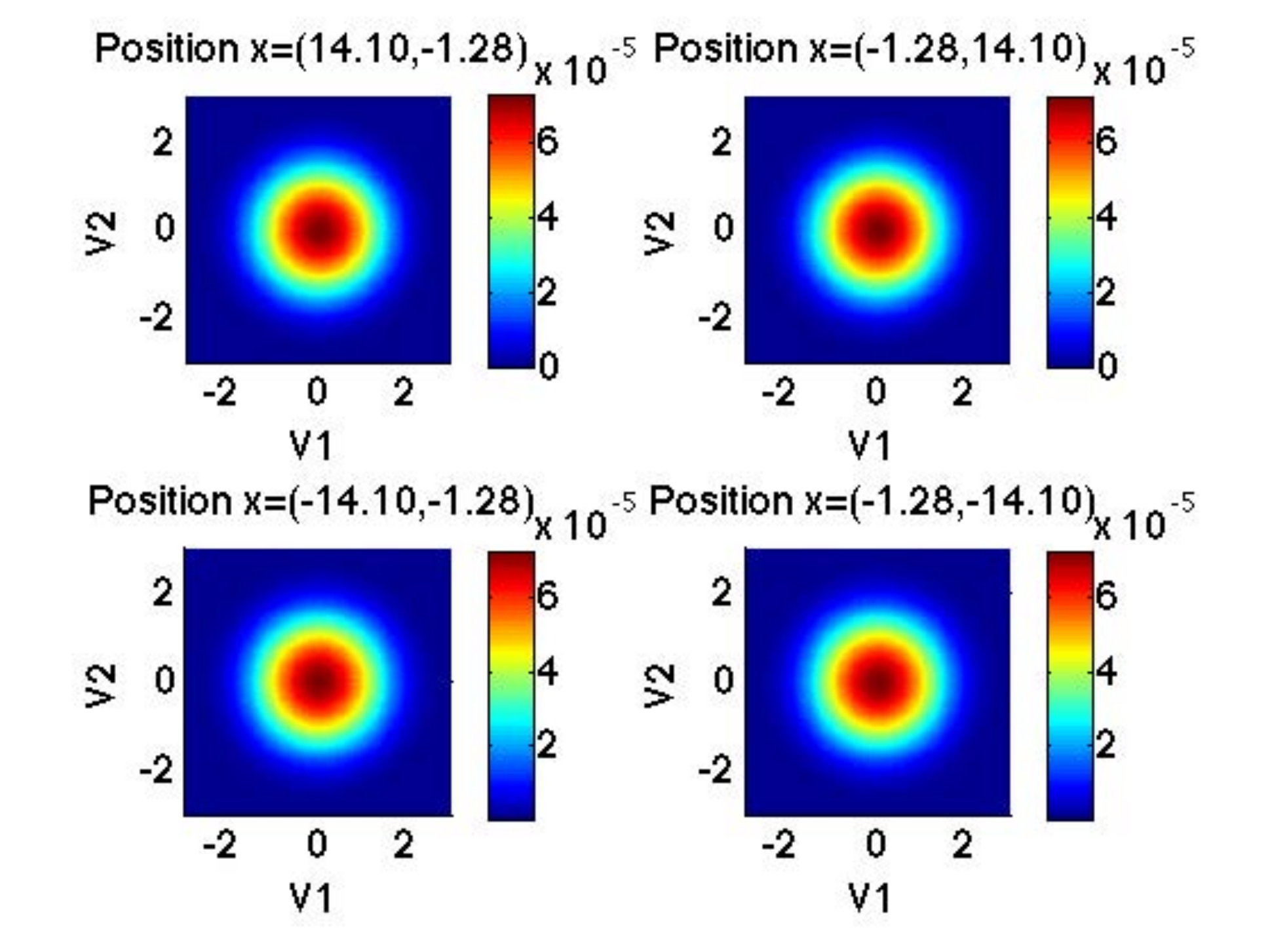}
\caption{Solutions of the mean field  equation: On the left, we have the spatial density $\int f\;dv$ and on the right, we have the velocity distribution at some fixed positions in the spatial domain. From top to bottom there are results for different values of the diffusion coefficient $A=0,0.123,0.9$.  The constants of the interaction potential are chosen according to \cite{hauptpaper}: $C_a=20,C_r=50,l_a=100,l_r=2$. \label{dmillfp}}
\end{figure}

The crucial point is, that starting from the single mill state for $A=0$, there is a transition from single to double mills for small values of $A$, before milling structures vanish as $A$ increases.
This can be observed for the solutions of the kinetic equation as well. Figure \ref{dmillfp} shows spatial density and velocity distributions at fixed points in space, for several values of $A$ for single mill initial conditions. The evaluation points in space are indicated above the respective velocity distribution. For $A=0$ one observes the milling solutions. In this case a $\delta$-type velocity distribution is observed for the fixed spatial points and a circular density distribution with zero density in the center, as before. However, $A=0.123$  yields a double mill type stationary state, as can be seen by comparing the result to the double mill for $A=0$ in figure \ref{doublemillkinmiccomp}. Further increase of $A$ yields the disappearance of milling structures.

\begin{figure}
\centering
\includegraphics[width=6.3cm]{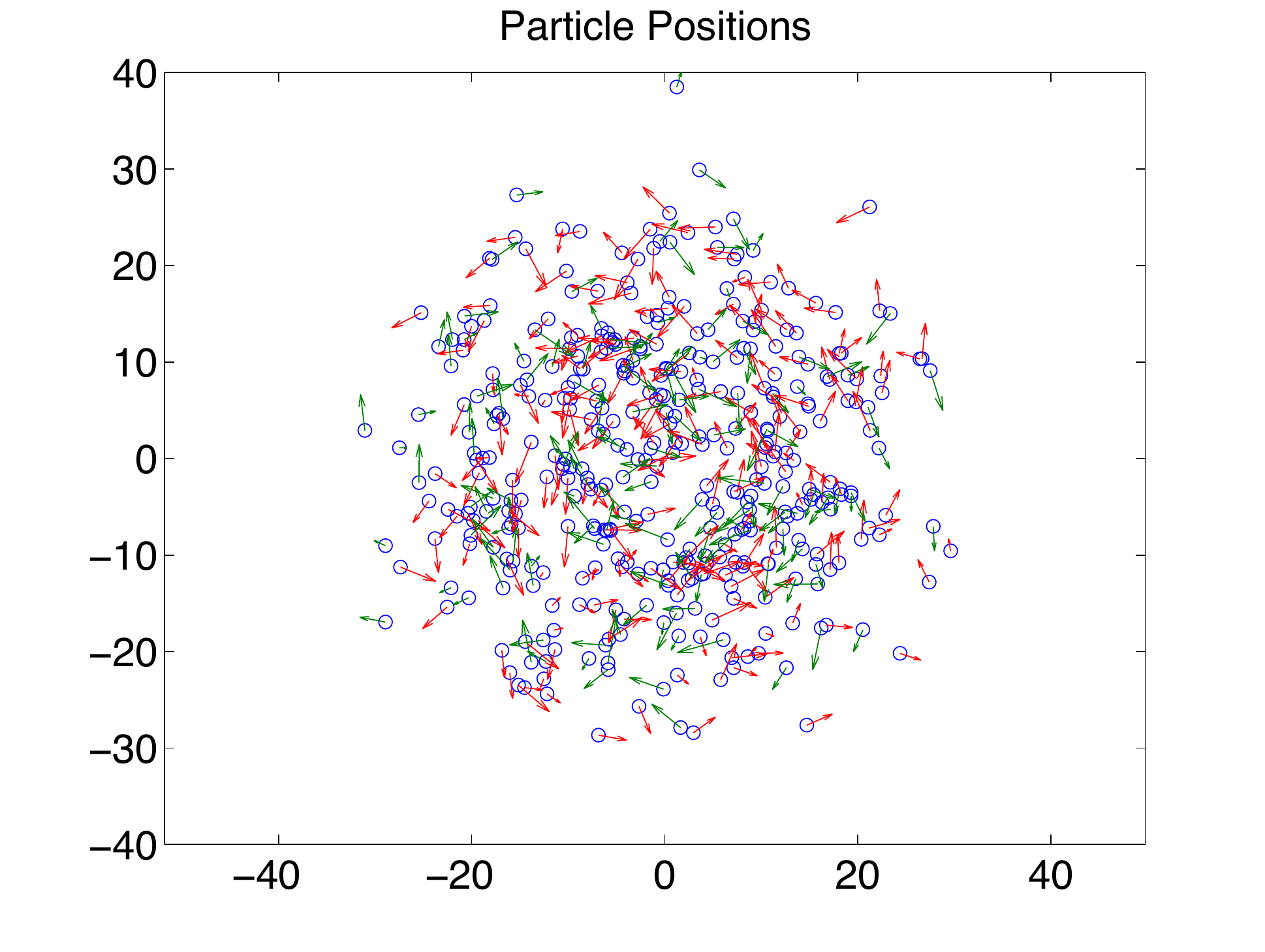}
\includegraphics[width=6.3cm]{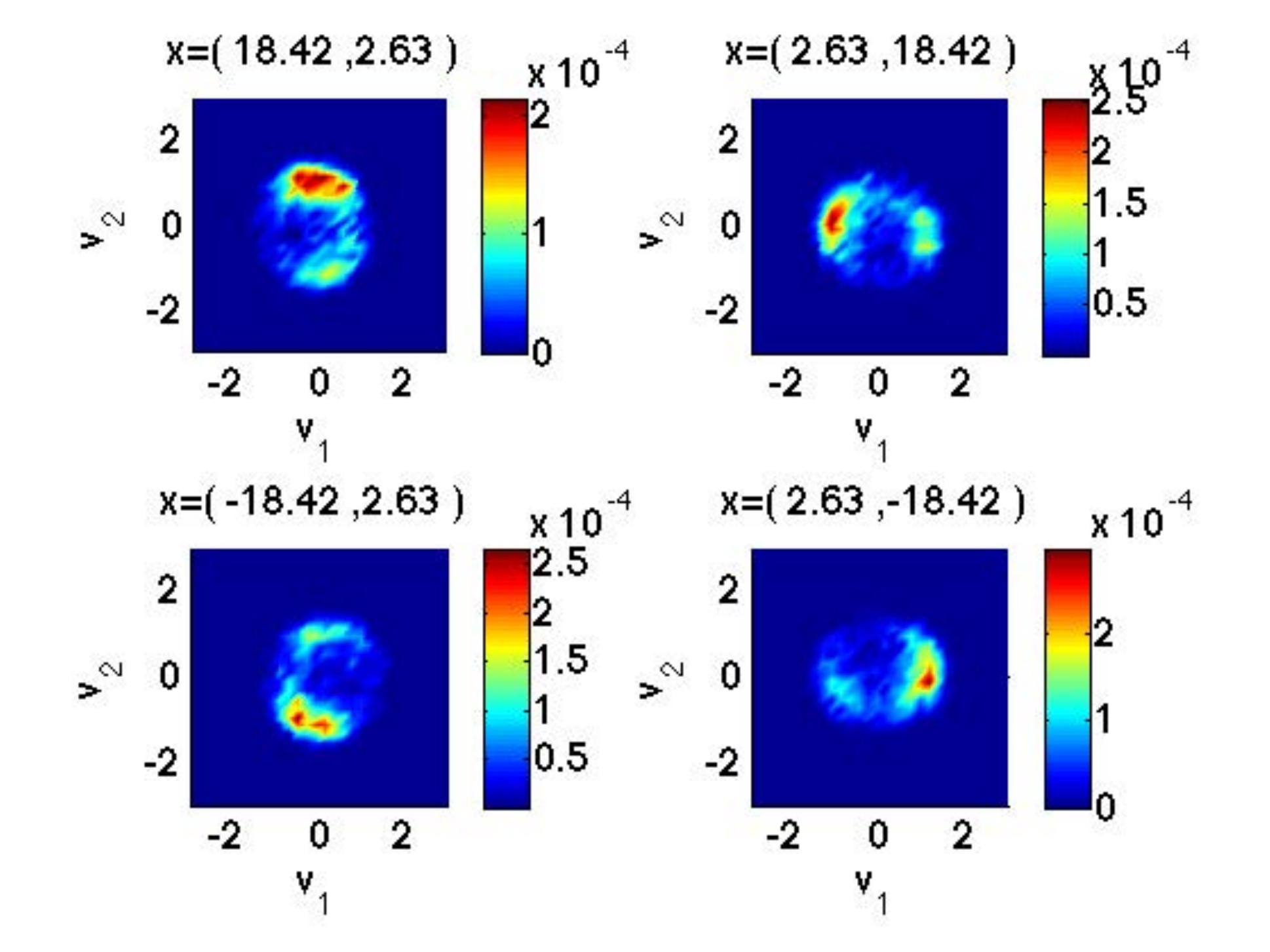}\\
\includegraphics[width=6.3cm]{kineticSpatialDensA0123.pdf}
\includegraphics[width=6.3cm]{kineticVelocitiesA0123.pdf}
\caption{Comparison of results from kinetic equation (top) and from microscopic equations (bottom). In each case, a single mill initial condition was used and the diffusion/noise parameter was set to $A=0.123$.}\label{a0123kinmiccomp}
\end{figure}

Moreover,  for $A=0.123$ we compare the microscopic and the kinetic solutions. We  generate again  a histogram from the particle data and give the comparison to the kinetic result in figure \ref{a0123kinmiccomp}.  Again, the solutions match, with a slightly higher diffusivity in  the kinetic result  due to the numerical scheme.

\begin{figure}
\centering
\includegraphics[width=6.3cm]{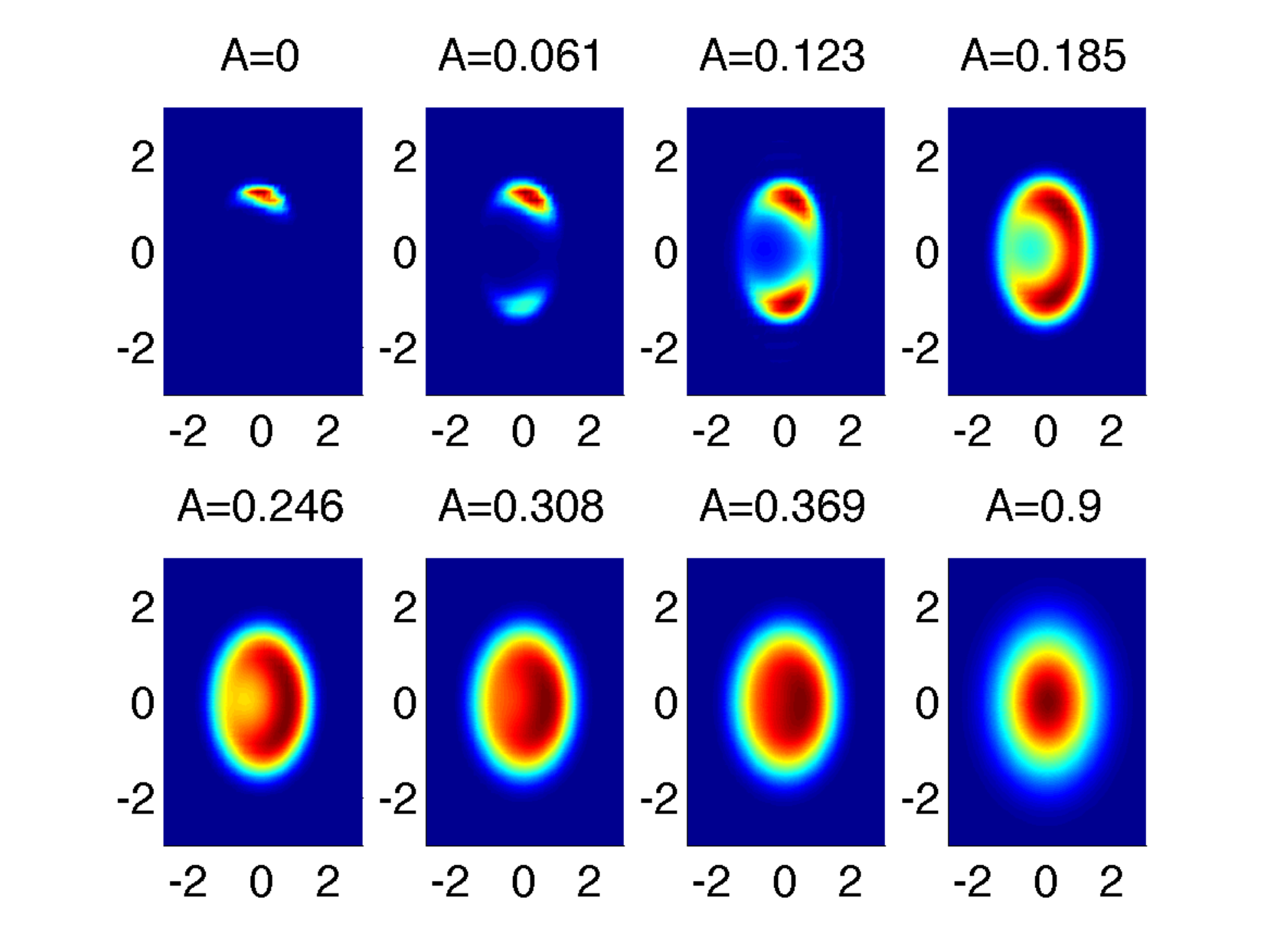}
\includegraphics[width=6.3cm]{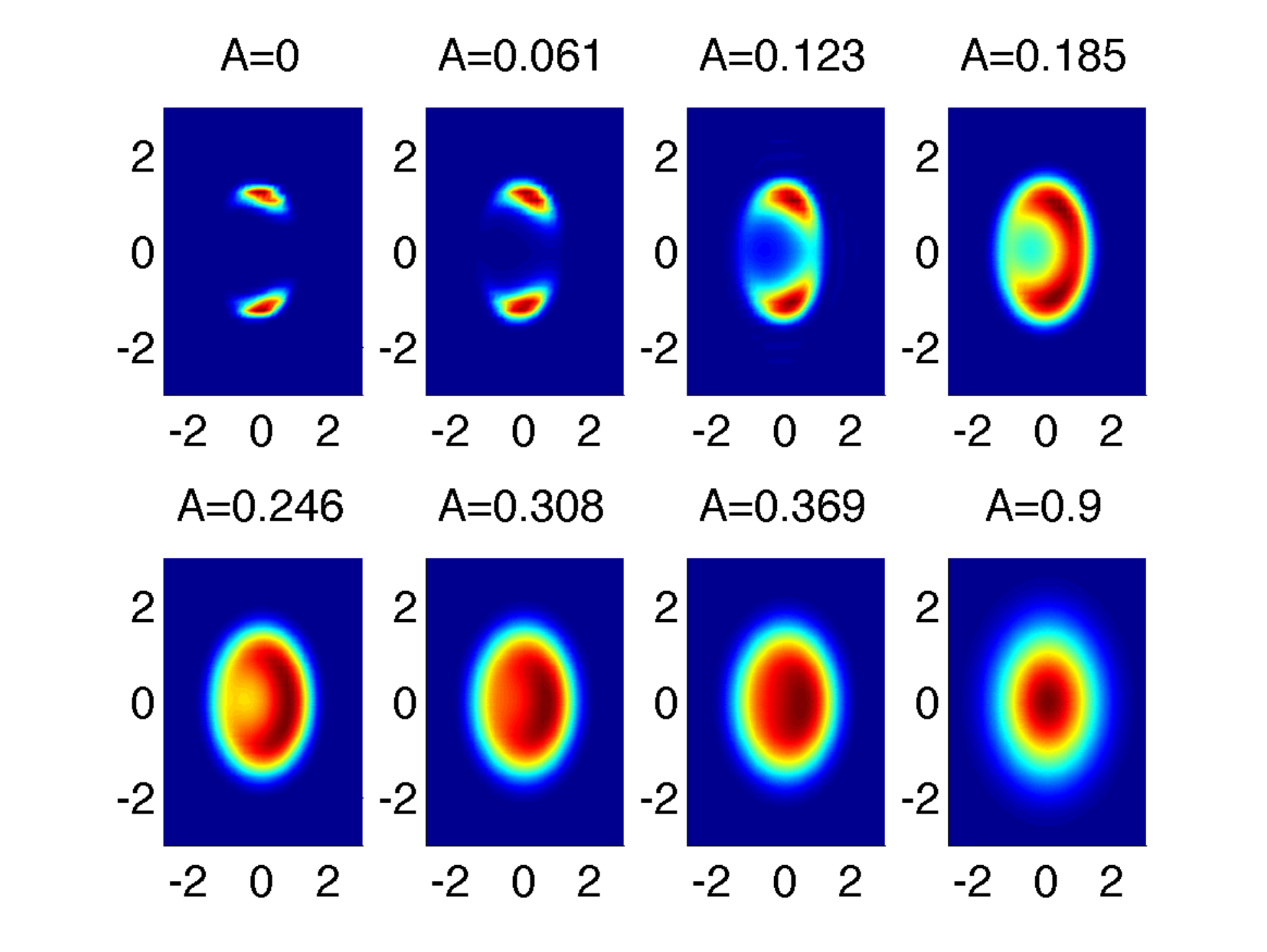}
\caption{Here, we look at $f$ at one of the fixed spatial points. On the left, we have a single mill, which fades away gradually with increasing diffusion coefficient $A$. On the right, we have the same situation starting from a double mill.\label{diverseA}}
\end{figure}

Additionally, we plot in Figure \ref{diverseA} the velocity distribution functions at one fixed spatial point for a wider range of values of $A$, with separate plots for single mill and double mill initial conditions. One can see, that the single mill state fades into a double mill and then into a normal distribution. The double mill however persists for small values of $A$.

\begin{figure}
\centering
\includegraphics[width=8cm]{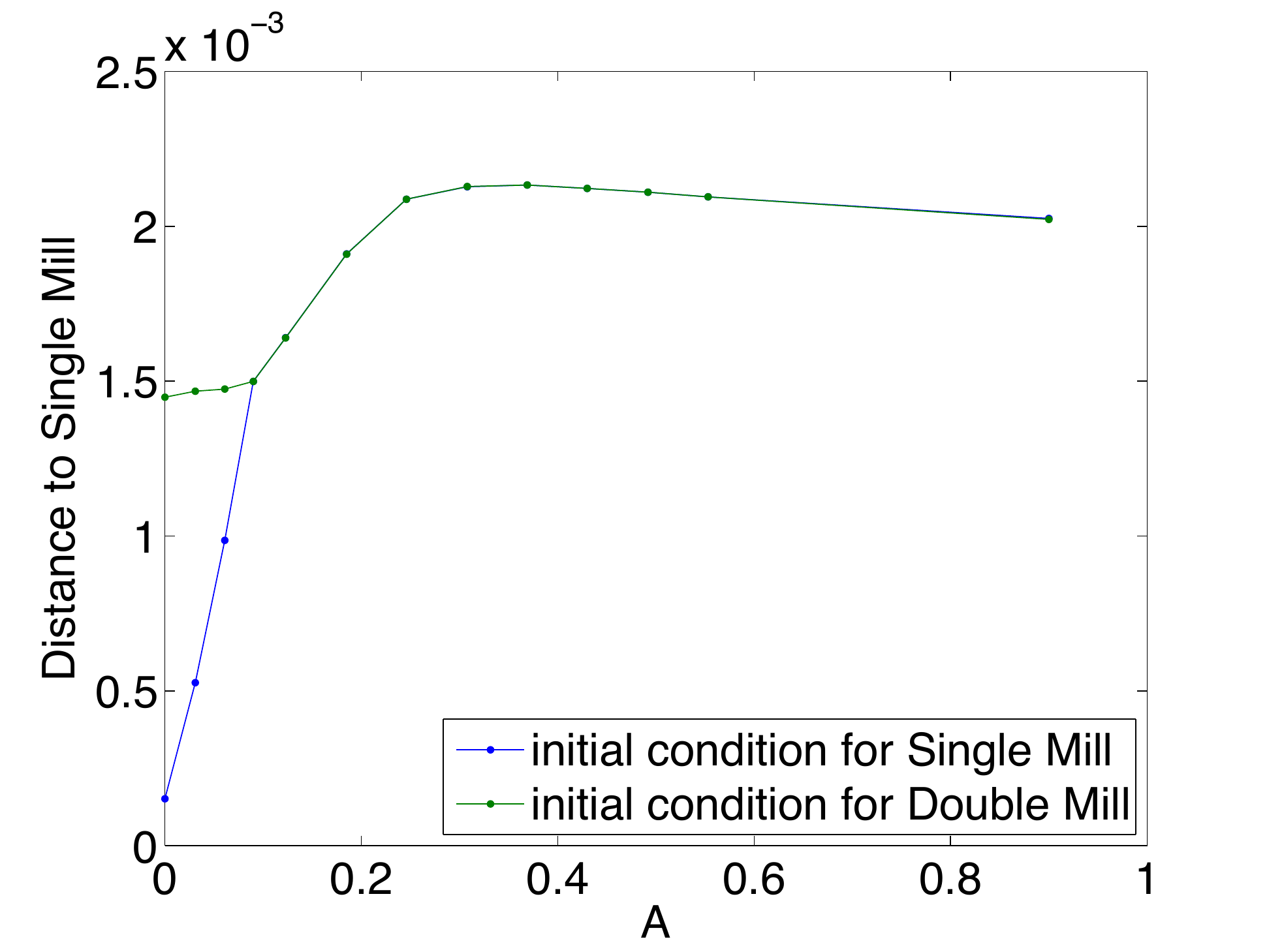}
\caption{Distance of the solution from the single  mill solution (\ref{distmeas}) for different values of $A$. The blue line is obtained for single mill initial condition, the green one for a double mill initial condition.
\label{milldist}}
\end{figure}

In order to examine this behaviour more closely, we computed for each numerical solution $f_A$ of the kinetic equation the distance to a single mill
\begin{align}
\max_{x \in \mbox{supp} \rho} \int \left| v - \sqrt{\frac{\alpha}{\beta}} \frac{x^\perp}{\vert x \vert} \right|^2 f_A (x,v) dv .
\label{distmeas}
\end{align}
Figure \ref{milldist} shows the values of this functional depending on the diffusion parameter $A$ for initial conditions yielding either single or double mills. We can observe the presumed behaviour, where the double mill persists for some range of $A$ and then fades away as the single mill does.

\begin{figure}
\centering
\hspace*{-2cm}\includegraphics[width=16.5cm]{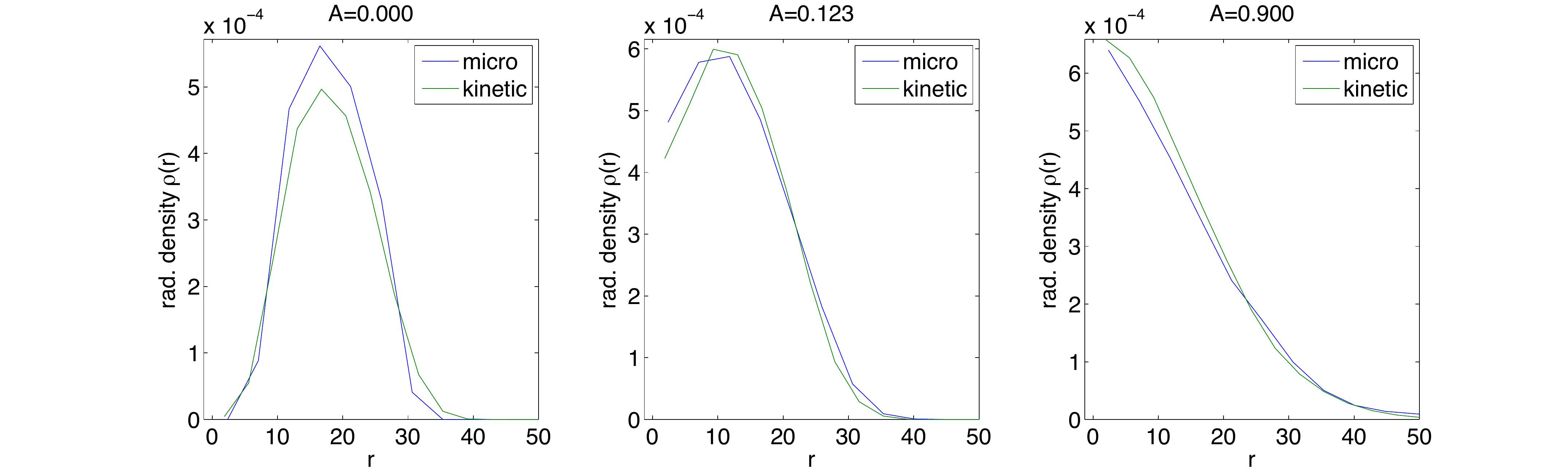}
\caption{Comparison of spatial densities from microscopic and kinetic equations. We show radial plots of them for different values of $A$: From left to right, we have $A=0.0$, $A=0.369$ and $A=0.9$.\label{radial}}
\end{figure}

\begin{figure}
\centering
\includegraphics[width=13cm]{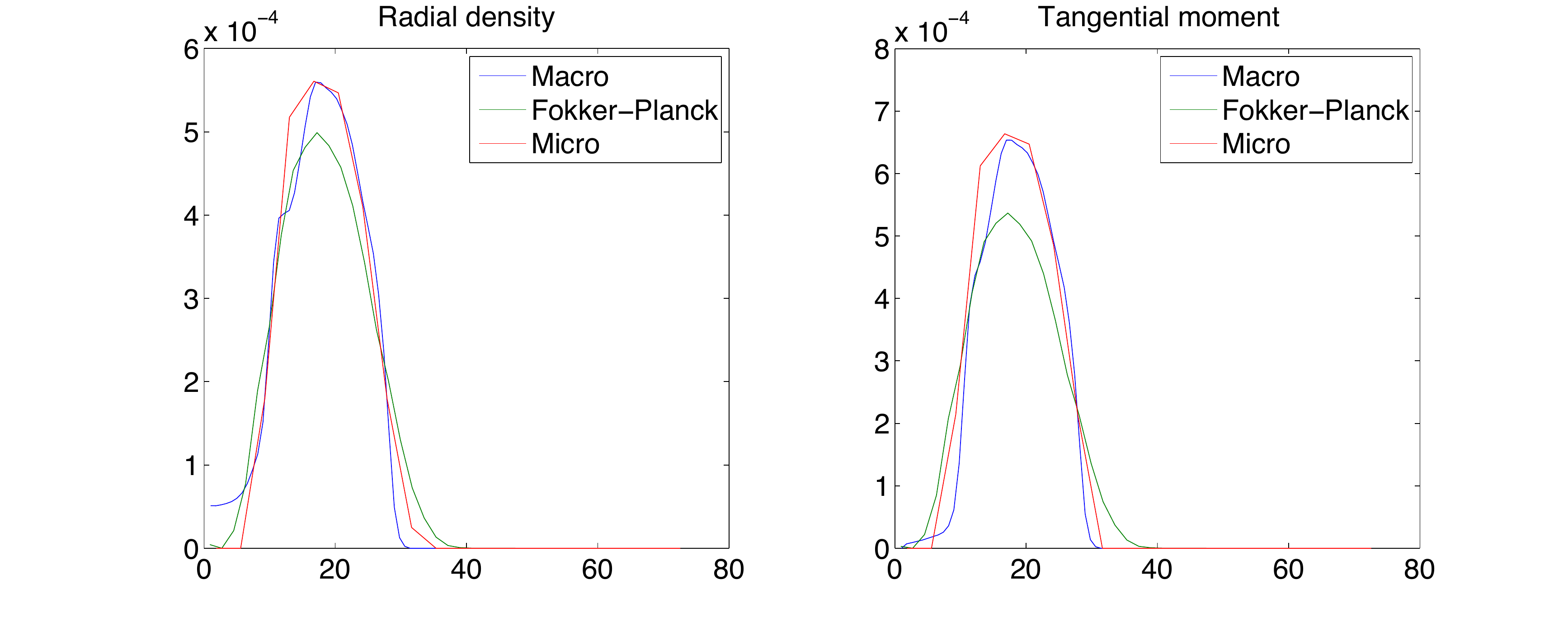}
\caption{On the left: Radial plot of the densities obtained from microscopic, mean field and macroscopic single mill solution for $A=0$. On the right: radial plot of the corresponding tangential moment. \label{micromacro}}
\end{figure}

\begin{figure}
\centering
\includegraphics[width=10cm]{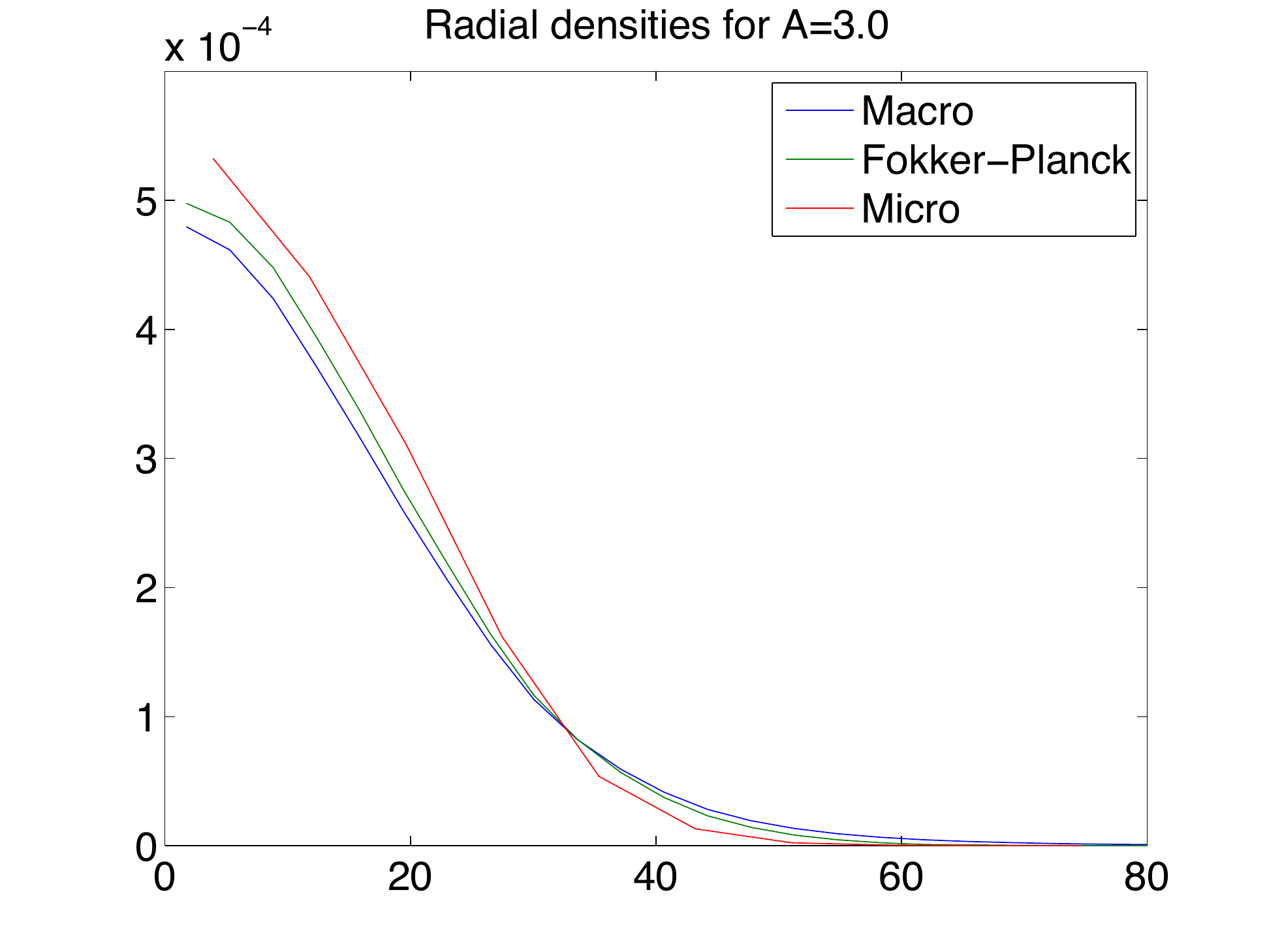}
\caption{Radial plot of the densities obtained from microscopic, mean field and macroscopic solution for  $A=3.0 $.\label{Abig}}
\end{figure}

In order to further validate our results for the kinetic equation, we have also compared the spatial densities of microscopic, kinetic and macroscopic equations. In figure \ref{radial} one can find radial plots of the spatial densities of the mean-field equation and the microscopic system for different values of $A$. For $A=0$ we have the above-mentioned circular distribution with zero density in the center, which fades away for $A=0.123$ until there is a normal distribution for $A=0.9$ and higher values of $A$. Furthermore, we compare our microscopic and kinetic solutions to the macroscopic single mill result computed in \cite{hauptpaper} for $A=0$ in figure \ref{micromacro}. We get good agreement of both the spatial density $\rho$ and the tangential moment $\rho u$. The mean-field solution is a bit more diffusive than the corresponding microscopic and macroscopic ones. We note again that the macroscopic solution is not able to reproduce  the double mill solution, since a mono-kinetic closure is used in section \ref{A=0}. For large $A$ the radial density  of the mean field solution is shown together with the corresponding macroscopic solution of  the diffusion problem from  section \ref{Alarge} and the microscopic solution in figure \ref{Abig}.

\subsection{Computational remarks}

All the computations were carried out on the cluster of the University of Kaiserslautern on an Intel Xeon E5 2670 with eight cores at 2.6 GHz each. The solver for the Fokker-Planck equation is written in Fortran and with OpenMP parallelisation. We note that for the mean field equation, we only need to compute  distances between points of the  grid to determine the influence of the interactions. Since the grid is fixed, these distances can be precomputed.   The methods discussed here can be  extended to three dimensions in  space and velocity domain, but then memory consumption and computing time will become  important issues. For example, the storage of the precomputed  pairwise distances between grid points in 3 dimensions takes a considerable amount of memory.  

\subsection*{Acknowledgements}
JAC acknowledges support from the MINECO Spanish projects MTM2011-27739-C04-02 (FEDER),
2009-SGR-345 from Ag\`encia de Gesti\'o d'Ajuts Universitaris i de Recerca-Generalitat de Catalunya, the Royal Society through a Wolfson Research Merit Award, and the Engineering and Physical Sciences Research Council (UK) grant number EP/K008404/1. The work of
the last two authors has been supported by the German research foundation, DFG grant KL
1105/17-1 and KL 1105/18-1.


\begin{thebibliography}{10}

\bibitem{ABCV}
Albi, G., Balagu\'e, D., Carrillo, J. A., von Brecht, J.: \emph{Stability Analysis of Flock and Mill Rings for Second Order Models in Swarming}. SIAM J. Appl. Math., {\bf 74} (2014), pp. 794--818. 

\bibitem{rome}
Ballerini, M., Cabibbo, N., Candelier, R., Cavagna, A., Cisbani,
E., Giardina, L., Lecomte, L., Orlandi, A., Parisi, G.,
Procaccini, A., Viale, M., Zdravkovic, V.: \emph{Interaction ruling
animal collective behavior depends on topological rather than
metric distance: evidence from a field study}. Proc Natl Acad Sci USA, \textbf{105} (2008) , pp. 1232-1237.

\bibitem{BTTYB} Barbaro, A., Taylor, K., Trethewey, P.F.,
Youseff, L., Birnir, B.: \emph{Discrete and continuous models of the
dynamics of pelagic fish: application to the capelin}. Mathematics
and Computers in Simulation, \textbf{79} (2009), pp. 3397--3414.

\bibitem{BEBSVPSS} Barbaro, A., Einarsson, B., Birnir, B., Sigurthsson, S.,
Valdimarsson, H., Palsson, O.K.,  Sveinbjornsson, S., Sigurthsson,
T.: \emph{Modelling and simulations of the migration of pelagic fish}.
ICES J. Mar. Sci., \textbf{66} (2009), pp. 826--838.

\bibitem{bigring}
Bertozzi, A. L., von Brecht, J., Sun, H., Kolokolnikov, T., Uminsky, D.: \emph{Ring Patterns and their Bifurcations in a Nonlocal Model of Biological Swarms}, to appear in Comm. Math. Sci., (2014).

\bibitem{BCC}
Bolley, F., Canizo, J.A., Carrillo, J.A.: \emph{Stochastic Mean-Field Limit: Non-Lipschitz Forces \& Swarming}. Math. Mod. Meth. Appl. Sci., {\bf 21} (2011), pp. 2179--2210. 

\bibitem{BDT}
Bonabeau, E., Dorigo, M., Theraulaz, G.: \emph{Swarm Intelligence: From
Natural to Artificial Systems}. Intelligence: From Natural to Artificial Systems (Oxford University Press, New York, 1999);

\bibitem{BGKMW07}
Bonilla, L., G\"otz, T., Klar, A., Marheineke, N., Wegener, R.:
\emph{Hydrodynamic limit for the Fokker-Planck equation of fiber lay--down models}.
SIAM Appl. Math., \textbf{68} (2007), pp. 648-655.

\bibitem{BH} Braun, W., Hepp, K.: \emph{The Vlasov Dynamics
and Its Fluctuations in the 1/N Limit of Interacting Classical
Particles}. Commun. Math. Phys., \textbf{56} (1977),  pp. 101-113.

\bibitem{camazine}
Camazine, S., Deneubourg, J.-L., Franks, N.R., Sneyd, J.,
Theraulaz, G., Bonabeau, E.: Self-Organization in Biological
Systems. Princeton University Press  (2003).

\bibitem{CCR}
Ca\~nizo, J.A.,  Carrillo, J.A., Rosado, J.: \emph{A well-posedness
theory in measures for some kinetic models of collective motion}.
Math. Mod. Meth. Appl. Sci., {\bf 21} (2011), pp. 515-539.

\bibitem{review2}
Ca\~nizo, J.A., Carrillo, J.A., Rosado, J.: \emph{Collective Behavior of
Animals: Swarming and Complex Patterns}. Arbor, {\bf 186} (2010), pp. 1035--1049.

\bibitem{CDP}
Carrillo, J.A., D'Orsogna, M.R., Panferov, V.: \emph{Double milling in
self-propelled swarms from kinetic theory}. Kinetic and Related
Models, \textbf{2} (2009), pp. 363-378.

\bibitem{CHM2}
Carrillo, J. A., Huang, Y., Martin, S.: \emph{Nonlinear stability of flock solutions in second-order swarming models}. Nonlinear Anal. Real World Appl., {\bf 17} (2014), pp. 332--343.

\bibitem{CHM}
Carrillo, J. A., Huang, Y., Martin, S.: \emph{Explicit flock solutions for quasi-morse potentials}.
to appear in European J. Appl. Math., (2014).

\bibitem{CFRT}
Carrillo, J.A., Fornasier, M., Rosado, J., Toscani, G.: \emph{Asymptotic
Flocking Dynamics for the kinetic Cucker-Smale model}. SIAM J. Math. Anal., {\bf 42} (2010), pp. 218--236.

\bibitem{review}
Carrillo, J.A., Fornasier, M., Toscani, G., Vecil, F.:  \emph{Particle,
Kinetic, and Hydrodynamic Models of Swarming}, in Naldi, G., Pareschi, L., Toscani, G. (eds.) Mathematical Modeling of Collective Behavior in Socio-Economic and Life Sciences, Series: Modelling and Simulation in Science and Technology, Birkhauser, (2010), pp. 297--336. 

\bibitem{hauptpaper}
Carrillo, J.A., Klar, A., Martin, S., Tiwari, S.: \textit{Self-propelled interacting particle systems with roosting force}. Math. Mod. Meth. Appl. Sci., {\bf 20} (2010), pp. 1533--1552.

\bibitem{Carrillo2013}
Carrillo, J.A., Martin, S., Panferov, V.: \emph{A new interaction potential for swarming models}.
Physica D, {\bf 260} (2013), pp. 112--126.

\bibitem{CV}
Carrillo, J. A., Vecil, F.: \emph{Nonoscillatory interpolation methods applied to {V}lasov-based models}. SIAM J. Sci. Comput., {\bf 29} (2007), pp. 1179--1206.

\bibitem{Splittingstrang}
Cheng, C., Knorr, G.: {\em The integration of the Vlasov equation in configuration space}, 
J.Comp.Phys., {\bf 22} (1976), pp. 330--351.

\bibitem{CDMBC}
Chuang, Y.L., D'Orsogna, M.R., Marthaler, D., Bertozzi, A.L.,
Chayes, L.: \emph{State transitions and the continuum limit for a 2D
interacting, self-propelled particle system}. Physica D,
\textbf{232} (2007), pp. 33-47.

\bibitem{couzin}
Couzin, I.D., Krause, J., Franks, N.R., Levin, S.A.: \emph{Effective
leadership and decision making in animal groups on the move}.
Nature, \textbf{433} (2005), pp. 513-516.

\bibitem{CKJRF} Couzin, I.D., Krause, J., James, R., Ruxton, G. and Franks, N.:
\emph{Collective memory and spatial sorting in animal groups}. Journal of
Theoretical Biology, \textbf{218}  (2002), pp. 1-11.

\bibitem{DFL}
Degond, P., Frouvelle, A., Liu, J.-G.: \emph{Macroscopic limits and phase transition in a system of self-propelled particles}. J. Nonlinear Sci., {\bf 23} (2013), pp. 427--456. 

\bibitem{DM1} Degond, P., Motsch, S.:
\emph{Continuum limit of self-driven particles with orientation
interaction}. Math. Models Methods Appl. Sci.,  \textbf{18} (2008), pp.
1193-1215.

\bibitem{dobru}
Dobrushin, R.: \emph{Vlasov equations}. Funct. Anal. Appl., \textbf{13} (1979),
pp. 115-123.

\bibitem{DCBC}
D'Orsogna, M.R., Chuang, Y.L., Bertozzi, A.L., Chayes, L.:
\emph{Self-propelled particles with soft-core interactions: patterns,
stability, and collapse}. Phys. Rev. Lett., \textbf{96} (2006), 104302

 \bibitem{douglas00}  
Douglas, J., Huang, C., Pereira, F.: \emph{The Modified Method of Characteristics with Adjusted Advection}. Numer. Math., {\bf83} (1999), pp. 353--369.

\bibitem{filbet}
Filbet, F., Sonnendr\"ucker, E., Bertrand, P.: {\em Conservative
numerical schemes for the Vlasov equation}. J. Comput. Phys., {\bf 172}
(2001), pp. 166--187.

\bibitem{GKMW07}
G\"otz, T., Klar, A., Marheineke, N., Wegener, R.:
\emph{A stochastic model and
associated Fokker--Planck equation for the fiber lay-down process
in nonwoven production processes}. SIAM Appl. Math, \textbf{67} (2007), pp. 1704-1717

\bibitem{herleitungallgemein}
Golse, F.: 
\textit{The mean field limit for the dynamics of large particle systems}.
Journ\'ees \'equations aux d\'eriv\'ees partielles, {\bf 9} (2003), pp. 1--47.

\bibitem{chate}
Gr\'egoire, G., Chat\'e, H.: \emph{Onset of collective and cohesive
motion}. Phy. Rev. Lett., \textbf{92} (2004), 025702

\bibitem{HL08}
Ha, S.-Y., Liu, J.-G.: \emph{A simple proof of the Cucker-Smale flocking
dynamics and mean-field limit}. Comm. Math. Sci., \textbf{7}, (2009), pp. 297-325

\bibitem{sem6}
Sonnendr{\"u}cker, E., Roche, J., Bertrand, P., Ghizzo, A.: \emph{The Semi-Lagrangian Method for the Numerical Resolution of the Vlasov Equation}. J. Comp. Phys. {\bf 149} (1999) pp. 201--220.
	
\bibitem{HT08}
Ha, S.-Y., Tadmor, E.: \emph{From particle to kinetic and hydrodynamic
descriptions of flocking}. Kinetic and Related Models, \textbf{1} (2008),
pp. 415-435.

\bibitem{HK} Hemelrijk, C.K. and Kunz, H.: \emph{Density distribution and
size sorting in fish schools}: an individual-based model.
Behavioral Ecology, \textbf{16} (2005), pp. 178-187.

\bibitem{HKMO09} Herty, M., Klar, A., Motsch, S., Olawsky, F.:
\emph{A smooth model for fiber lay-down processes and its diffusion approximations}. Kinetic and Related Models \textbf{2} (2009), pp. 480-502.

\bibitem{HCH} Hildenbrandt, H, Carere, C., and Hemelrijk, C. K.:
\emph{Self-organised complex aerial displays of thousands of starlings:
a model}. Behav. Ecol., {\bf 21} (2010), pp. 1349--1359.

\bibitem{HW} Huth, A. and Wissel, C.: \emph{The Simulation of the Movement of
Fish Schools}. Journal of Theoretical Biology, \textbf{152} (1992), pp. 365-385

\bibitem{KRS07} Klar, A., Reutersw\"ard, P., Sea\"{\i}d, M.:
\emph{A semi-Lagrangian method for the Fokker--Planck equation of fiber dynamics}.
J. Sci. Comput., \textbf{38} (2009), pp. 349-367

\bibitem{3dfadenmodell}
Klar, A., Maringer, J., Wegener, R.: \textit{A 3D model for fiber lay-down in nonwoven production processes}. Math. Mod. Meth. in Appl. Sci, {\bf 22(9)} (2012).

\bibitem{semilagrange}
Klar, A., Reutersw\"ard, P., Sea\"{\i}d, M.:
{\em A semi-Lagrangian method for a Fokker-Planck equation describing fiber dynamics}.
J.Sci.Comp., {\bf 38} (2009), pp. 349--367.

\bibitem{kolo}
Kolokolnikov, T., Carrillo, J. A., Bertozzi, A. L., Fetecau, R., Lewis, M.:
{\em Emergent behaviour in multi-particle systems with non-local interactions}.
{Phys. D}, {\bf 260} (2013), pp. {1--4}.

\bibitem{KH} Kunz, H. and Hemelrijk, C. K.: \emph{Artificial fish schools:
collective effects of school size, body size, and body form}. Artificial Life, \textbf{3} (2003), pp. 237-253.

\bibitem{leveque}
Leveque, R.J.: 
\textit{Finite Volume Methods for Hyperbolic Problems}.
Cambridge University Press (2002).

\bibitem{LR}
Levine, H., Rappel, W.J., Cohen, I.: \emph{Self-organization in systems
of self-propelled particles}. Phys. Rev. E, \textbf{63} (2000), 017101.

\bibitem{LLE}
Li, Y.X., Lukeman, R., Edelstein-Keshet, L.: \emph{Minimal mechanisms
for school formation in self-propelled particles}. Physica D,
\textbf{237} (2008), pp. 699-720.

\bibitem{LLE2}
Li, Y.X., Lukeman, R., Edelstein-Keshet, L.: \emph{A conceptual model
for milling formations in biological aggregates}. Bull Math Biol.,
\textbf{71} (2008), pp. 352-382.

\bibitem{mogilner}
Mogilner, A., Edelstein-Keshet, L., Bent, L., Spiros, A.: \emph{Mutual
interactions, potentials, and individual distance in a social
aggregation}. J. Math. Biol., \textbf{47} (2003), pp. 353-389.

\bibitem{MT}
Motsch, S., Tadmor, E.:
A new model for self-organized dynamics and its flocking behavior. 
J. Stat. Phys., {\bf 144} (2011), pp. 923--947. 

\bibitem{neunzert}
Neunzert, H.: \emph{The Vlasov equation as a limit of Hamiltonian
classical mechanical systems of interacting particles}. Trans.
Fluid Dynamics, \textbf{18} (1977), pp. 663-678.


\bibitem{parrish}
Parrish, J., Edelstein-Keshet, L.: \emph{Complexity, pattern, and
evolutionary trade-offs in animal aggregation}. Science, {\bf 294} (1999),
pp. 99-101.

\bibitem{quarteroni}
{Quarteroni, A., Valli, A.}: 
{Numerical Approximation of Partial Differential Equations}.
Springer, Second Edition (1997).

\bibitem{qiu}
{Qiu, J.-M., Christlieb, A.}:
{\em A Conservative high order semi-Lagrangian WENO method for the Vlasov Equation}.
Journal of Computational Physics, {\bf 229} (2010), pp. 1130--1149.

\bibitem{numerikfaden}
{Roth, A., Zharovsky, E., Klar, A., Simeon,B.}:
\textit{ A Semi-Lagrangian Method for 3-D Fokker Planck Equations for Stochastic Dynamical Systems on the Sphere}.
to appear in J. Sci. Comput., (2014).

\bibitem{spohn2}
Spohn, H.: \emph{Large scale dynamics of interacting particles}. Texts
and Monographs in Physics, Springer (1991).

\bibitem{herleitungallgemein2}
{Spohn, H.}:
\textit{Kinetic equations from Hamiltonian dynamics: Markovian limits}.
Rev. Modern Phys., {\bf 52} (1980), pp. 569--615.


\bibitem{toro}
{Toro, E.}:
\textit{Riemann Solvers and Numerical Methods for Fluid Dynamics}.
Springer, Third Edition (2009).

\bibitem{vanleer}
{Van Leer, B.}:
\textit{Towards the ultimate conservative difference scheme II. Monotonicity and conservation combined in a second order scheme}.
J. Comp. Phys., {\bf 14} (1974), pp. 361--370.

\bibitem{VLR}
Vecil, F., Laffitte, P., Rosado, J.: \emph{A numerical study of attraction/repulsion collective behavior models: 3D particle analyses and 1D kinetic simulations}. Physica D., {\bf 260} (2013), pp. 127--144.

\bibitem{VCBCS}
Vicsek, T., Czirók, A., Ben-Jacob, E., Cohen, I., Shochet, O.: \emph{Novel type of phase transition in a system of self-driven particles}. Phys. Rev. Lett., {\bf 75} (1995), pp. 1226--1229. 

\bibitem{evgeniy}
{Zharovsky, E., Simeon, B.}: 
\textit{A space-time adaptive approach to orientation dynamics in particle laden flow}. 
Procedia Computer Science, {\bf 1} (2010), pp. 791--799.

\bibitem{evgeniyhauptpaper}
{Zharovski, E., Moosaie, A., LeDuc, A., Manhart, M., Simeon, B.}:
\textit{On the numerical solution of a convection-diffusion equation for particle orientation dynamics on geodesic grid}.
IMACS Journal Applied Numerical Mathematics, {\bf 62} (2012), pp. 1554--1566.

\end{thebibliography}
\end{document}